\documentclass[12pt]{amsart}
\usepackage{amssymb,  amsfonts, amsmath, xy,latexsym, bbm}
\usepackage{mathrsfs}

\usepackage[margin=2cm]{geometry}
\usepackage[pdftex,colorlinks,hypertexnames=false,linkcolor=black,urlcolor=black,citecolor=black]{hyperref}

\numberwithin{equation}{section}

\newtheorem{Theorem}{Theorem}[section]
\newtheorem*{Theorem*}{Theorem}
\newtheorem*{Corollary*}{Corollary}

\newtheorem{Lemma}[Theorem]{Lemma}
\newtheorem{Proposition}[Theorem]{Proposition}
\newtheorem{Corollary}[Theorem]{Corollary}

\theoremstyle{definition}

\theoremstyle{remark}

\newtheorem*{Remark*}{Remark}

\newcommand{\C}{\mathbb{C}}

\newcommand{\Z}{\mathbb{Z}}
\newcommand{\Q}{\mathbb{Q}}
\renewcommand{\k}{\mathbbm{k}}
\renewcommand{\Z}{\mathbb{Z}}
\newcommand{\A}{\mathbb{A}}

\newcommand{\g}{\mathfrak{g}}

\newcommand{\gl}{\mathfrak{gl}}

\renewcommand{\t}{\mathfrak{t}}
\renewcommand{\l}{\mathfrak{l}}
\renewcommand{\sl}{\mathfrak{sl}}

\newcommand{\m}{\mathfrak{m}}

\newcommand{\D}{\mathcal{D}}
\renewcommand{\O}{\mathcal{O}}

\newcommand{\R}{\mathcal{R}}

\newcommand{\V}{\mathcal{V}}

\renewcommand{\L}{{\mathcal{L}}}

\newcommand{\ad}{\operatorname{ad}}

\newcommand{\Lie}{\operatorname{Lie}}
\newcommand{\Ker}{\operatorname{Ker}}

\newcommand{\SL}{\operatorname{SL}}
\newcommand{\GL}{\operatorname{GL}}

\newcommand{\Aut}{\operatorname{Aut}}

\newcommand{\End}{\operatorname{End}}

\newcommand{\Der}{\operatorname{Der}}
\renewcommand{\End}{\operatorname{End}}

\newcommand{\Char}{\operatorname{char}}

\newcommand{\del}{{\partial}}

\begin{document}
\title{The Harish-Chandra centre and the Azumaya locus of the second Witt--Jacobson Lie algebra}
\author{Chao Ma and Alexander Premet}
\thanks{\nonumber{\it Mathematics Subject Classification} (2020 {\it revision}).
Primary 17B35, 17B50. Secondary 17B66, 16H05.}
\address{Independent Researcher, London, UK}
\email{raymond.ma@me.com} 
\address{Department of Mathematics, The University of Manchester, Oxford Road, M13 9PL, UK} 
\email{alexander.premet@googlemail.com}
\maketitle
\begin{abstract}
\noindent 
Let $\O=\k[x,y]/(x^p,y^p)$ be the truncated polynomial ring over an algebraically closed field of characteristic $p>3$ and $\L=\Der\O$, the second Witt--Jacobson Lie algebra. 
Let $U$ be the universal enveloping algebra of $\L$. It is well known that $\L$ has index $2$ and the KW1 conjecture holds for $\L$, that is the maximal dimension of simple $\L$-modules equals $p^{(\dim\L-{\rm ind}\,\L)/2}=p^{p^2-1}$. It is also known that the field of fractions
of the centre $Z(U)$ is a separable extension of degree $p^2$ of the field of fractions of the $p$-centre $Z_p(\L)$. In this paper we describe the invariant ring
$S(\L)^\L$ and show that despite the lack of a symmetrisation map any homogeneous element  of $S(\L)^\L$ can be lifted to an element of $Z(U)$. We give an explicit presentation of $Z(U)$ and determine its singularities. As a consequence, we describe the Azumaya locus of $Z(U)$.
\end{abstract}

\section{Introduction and preliminaries}
\subsection{}
\label{ss:1.1}
Let $\k$ be an algebraically closed field of characteristic $p>3$ and $\O=\k[x,y]/(x^p,y^p)$.
The standard grading of the truncated polynomial ring $\O$ is obtained by giving the generators $x$ and $y$ degree $1$. The unique maximal ideal $\m$ of $\O$ is spanned by all monomials $x^iy^j$ with $i+j\ge 1$. Any element of $\O$ can be written as $f=f(0)+f'$ where $f(0)\in\k$ and $f'\in\m$.
The algebraic group $G=\Aut\O$ preserves $\m$. It is well known  (and easy to see) that $G=G_0\cdot R_u(G)$ where $G_0\cong \GL(2, \k)$ is the normaliser of  $\k x\oplus \k y$ and the unipotent radical $R_u(G)$ consists of all automorphisms $x\mapsto f(x,y), y\mapsto g(x,y)$ whose linear part is the identity.
It follows that $G$ is a connected group and $\Lie(G)\subset \Der \O$ preserves $\m$. Any element of $G$ can be presented uniquely as $g=g_0\cdot g_1$ with $g_0\in G_0$ and $g_1\in R_u(G)$. If $\varkappa\colon G\to \k^\times$ is a rational character of $G$ then $\ker\varkappa$ must contain the unipotent group $R_u(G)$ and the derived subgroup $(G_0,G_0)\cong \SL(2,\k)$ of $G_0$. In particular,  the central involution ${\rm diag}(-1,-1)$ of $(G_0,G_0)$ lies in $\ker \varkappa$. From this it is immediate that 
$\varkappa(g)=\det(dg_0)^n$ for all $g\in G$. Here $n=n(\varkappa)\in\Z$ and $\det(dg_0)$ is the determinant of the Jacobian matrix of the linear part of $g$.
\subsection{}
\label{ss:1.2}
The derivation algebra $\L:=\Der\O$ is a free $\O$-module of rank two generated by the partial derivatives $\del_x$ and $\del_y$. It is a restricted simple Lie algebra of Cartan type $W(2)$ with the $p$-mapping given by $D^{[p]}=D^p$ for all $D\in \L$. We write $\L_{(0)}$ for the standard maximal subalgebra $\m\del_x\oplus\m\del_y$ of $\L$. Jacobson proved in \cite{Jac} that under our assumption on $p$ 
each  automorphism of $\L$ is induced by a unique automorphism of $\O$, so that $\Aut\L\cong G$. This yields $\Lie(\Aut\L)=\ad\L_{(0)}$ and shows that $\Aut\L$ preserves $\L_{(0)}$ (in characteristic $2$ this is no longer true). The standard grading of $\O$ gives rise to that of $\L$, and we have $\deg(x^iy^j\del_x)=\deg(x^iy^j\del_y)=i+j-1$. Then $$\L=\oplus_{i=-1}^{2p-3}\,\L_i\ \ \ \ \ \mbox{and}\ \ \ \ \L^*=\oplus_{i=-1}^{2p-3}\,\L_i^*$$ where we embed the dual space $\L_i^*$ into $\L^*$ by letting $\oplus_{j\ne i}\,\L_j\subset \ker\varphi$ for all $\varphi\in\L_i^*$. Let $\{\del_x^*,\del_y^*\}$ be the basis of $\L_{-1}^*$ dual to the ordered basis $\{\del_x,\del_y\}$ of $\L_{-1}$.

Let $T_0$ be the $2$-dimensional maximal torus of $G$ generated  by the $1$-parameter subgroups $T_x\colon\,x\mapsto tx, y\mapsto y$ and  $T_y\colon\,x\mapsto x, y\mapsto ty$ with $t\in\k^\times$. We include $T_0$ into the Borel subgroup $B_0$ of $G_0$ with Lie algebra $\t_0\oplus \k x\del_y$ where $\t_0=\Lie(T_0)=\k x\del_x\oplus \k y\del_y$. Note that $\L_{-1}^*=\k \del_x^*\oplus\k \del_y^*$ and  $\L_{2p-3}=\k x^{p-1}y^{p-1}\del_x\oplus \k x^{p-1}y^{p-1}\del_y$ are irreducible $\L_0$-modules with highest weight vectors $\del_x^*$ and $x^{p-1}y^{p-1}\del_y$ with respect to $B_0$. The corresponding $\t_0$-weights are $(1,0)$ and $(-1,-2)$. 
\subsection{}
\label{ss:1.3}
Assume in this subsection that $\k$ is algebraically closed of characteristic $p>0$ and let $\g$ be an $n$-dimensional restricted Lie algebra over $\k$ with $p$-mapping $\g\ni x\to x^{[p]}\in\g$. The centre $Z$ of the universal enveloping algebra $U(\g)$ contains a large subalgebra generated by all $x^p-x^{[p]}$ with $x\in\g$. It is called the 
$p$-{\it centre} of $\g$ and denoted $Z_p(\g)$. The PBW theorem implies that if $\{e_i\,|\,\, 1\le i\le n\}$ is a $\k$-basis of $\g$ then $Z_p(\g)$ is a polynomial algebra  in $e_i^p-e_i^{[p]}$ and $U(\g)$ is a free $Z_p(\g)$-module of rank $p^n$. As a consequence, $U(\g)$ is a PI-domain and the maximal dimension $M(\g)=\mbox{PI-deg}(U(\g))$ of simple $U(\g)$-modules divides
$p^n$. 

To any simple $\g$-module $V$ one can attach a linear
function $\chi\in \g^*$, called the $p$-{\it character} of $V$, such that for any $x\in \g$ the element $x^p-x^{[p]}\in U(\g)$ acts on $V$ as
multiplication by $\chi(x)^p$. Conversely, any linear function on $\g$ serves as a $p$-character of a simple $\g$-module. This is immediate from the fact that for any $\chi\in\g^*$ the $\g$-modules with $p$-character $\chi$ are precisely the modules over the {\it reduced enveloping algebra} $$U_\chi(\g)\,:=\,U(\g)\otimes_{Z_p(\g)}\k_\chi\,\cong\,U(\g)/I_\chi,$$ where $\k_\chi=\k 1_\chi$ is the $1$-dimensional $Z_p(\g)$-module with $(x^p-x^{[p]})\cdot 1_\chi=\chi(x)^p1_\chi$ and $I_\chi$ is the two-sided ideal of $U(\g)$ generated by all $x^p-x^{[p]}-\chi(x)^p$ with $x\in \g$.

Let $Q(\g)={\rm Frac}\, Z$ and $Q_p(\g)={\rm Frac}\,Z_p(\g)$. 
It follows from results of Zassenhaus \cite{Zass} that $$[Q(\g):Q_p(\g)]=p^r \ \ \  \mbox{and}\ \ \  M(\g)=p^{(n-r)/2}$$ for some $r\in\Z_{\ge 0}$.
The first Kac--Weisfeiler conjecture \cite[1.2]{KW71} (abbreviated as KW1 following \cite{Kac}) predicts that $r={\rm ind}\,\g$, the {\it index} of $\g$. The latter is defined as 
 the minimal dimension of the coadjoint stabilisers $\g_\chi:=\{x\in\g\,\vert\,\, \chi([x,\g])=0\}$ with $\chi\in\g^*$. By \cite[Theorem~2]{Mil} and \cite[5.4]{PS99}, the inequality $r\le {\rm ind}\,\g$ holds in all cases.
 
 Since $U(\g)$ is a PI-domain, ${\rm Frac}\,U(\g)\,\cong\,U(\g)\otimes_{Z}\,Q(\g)$ is a central division algebra  over $Q(\g)$ of degree $M(\g)$ and there is a nonempty Zariski open subset $\Omega$ of ${\rm Specm}\,Z$ with the property that for $\m_\xi\in {\rm Specm}\,Z$ the algebra $U(\g)\otimes_{Z}\,\k_\xi$ is isomorphic to ${\rm Mat}_{M(\g)}(\k)$ if and only if $\xi\in \Omega$ (here $\xi\colon Z\to \k$ is the algebra homomorphism with kernel $\m_\xi$). The open set $\Omega$ is referred to as the {\it Azumaya locus} of $Z$. 
 
 Since $U(\g)$ is a finitely generated $Z_p(\g)$-module, so is $Z$. The inclusion $Z_p(\g)\hookrightarrow Z$ gives rise to a finite morphism $\Psi\colon\,{\rm Specm}\,Z\to {\rm Specm}\,Z_p(\g) \cong \g^*$ (suppressing the Frobenius twist is harmless at this point). As all fibres of $\Psi$ are finite and $\g^*$ is normal, $\Psi$ is an open morphism; see \cite[Ch.~AG, 18.4]{Bor}. Hence $\Psi(\Omega)$ is a nonempty Zariski open subset of $\g^*$.

 In some important cases we have the equality $\Psi(\Omega)=\g^*$. 
 This happens for instance when $\g=\Lie(G)$ where $G$ is a semisimple simply connected algebraic $\k$-group. Indeed, in this case analogues of the restricted Steinberg $\g$-module can be defined for all $p$-characters of $\g^*$ and the $U_\chi(\g)$-modules thus obtained are always irreducible of dimension $M(\g)=p^{(\dim G-{\rm rk}\, G)/2}$. This was first observed in \cite[Sect.~4]{FP90} by using deformation arguments (when $G$ is simply connected there is no need to impose any assumptions on $p$). However, when $\g=W(n)$ and $p>3$, the map $\Psi$ always sends $\Omega$ onto a {\it proper} open subset of $\g^*$. For $n=1$, this set is known thanks to  Chang's theorem \cite{Chang}. 
 Describing $\Psi(\Omega)$ for $W(n)$ with 
 $n\ge 2$ is an intriguing open problem. 
 
 For now, KW1 remains unproven in full generality and despite a plethora of objects to play with there are no known counterexamples (for restricted Lie algebras).
As mentioned in \cite[4.1]{PS99}, a solid ground for testing KW1 is provided by the Lie algebra $W(3)=\Der\, \k[x,y,z]/(x^p,y^p,z^p)$.
On the positive side, it is known that if there is $\chi\in\g^*$ such that $\g_\chi$ consists of $[p]$-semisimple elements, then necessarily $\g_\chi$ is a torus of dimension $r={\rm ind}\, \g$ and the field extension $Q(\g)/Q_p(\g)$ is  separable. Moreover, in this case there exists a nonempty Zariski open subset $\V$ of $\g^*$ such that for every $\psi\in \V$ the algebra $U_\psi(\g)$ is semisimple with $p^{{\rm ind}\,\g}$ pairwise non-isomorphic simple modules, each of dimension $M(\g)=p^{(n-{\rm ind}\,\g)/2}$; see \cite[Proposition~4.2(3) and Theorem~4.4]{PS99}.  In particular, KW1 holds in this case. It is conjectured in \cite[4.4]{PS99}
that the converse should be true as well, i.e. $U_\psi(\g)$ is semisimple for some $\psi\in\g^*$ only if $\g^*$ contains linear functions with toral stabilisers. 

We mention for completeness that Topley \cite{Top} gave an example of a finite dimensional completely solvable Lie algebra $\g$
for which $M(\g)>p^{(\dim \g-{\rm ind}\,\g)/2}$. However, in his example $\g$ does not admit a $p$-operation and one can show that KW1 still holds for a suitable $p$-envelope of $\g$
(as discussed in \cite[Chapter~1]{Str}). Although index is defined for any finite dimensional abstract Lie algebra, the notion of a $p$-centre has to be heavily modified in order to work in this setting; see \cite{Zass} for details.
\subsection{}
\label{ss:1.4} Let us now return to the Lie algebra $\L=W(2)$. The following facts were established in \cite{PS99}:

\begin{itemize}
	\item[1.] KW1 holds for $\L$ and $Q(\L)$ is a separable extension of $Q_p(\L)$ with  
	$[Q(\L)\colon Q_p(\L)]=p^2$;
	\item[2.] the Azumaya locus $\Omega$ and the smooth locus of the centre of $U(\L)$ coincide;
	\item[3.] there is $\eta\in \Psi(\Omega)$ such that $\L_\eta=\k\Delta \oplus \k\Delta^{[p]}$ where  $\Delta=\del_x-x^{p-1}\del_y$ is a regular nilpotent element of $\L$ (see \cite{Pr92, Pr15});
	\item[4.] the cocharacter $\mu\in X_*(T_0)$ 
	such that $\mu(t)\cdot x=tx$ and $\mu(t)\cdot y=t^p y$ has the property that $\mu(t)\cdot \eta=t^{{-(p^2-p-1})
	}\eta$ for all $t\in\k^\times$; 
	 \item[5.] the algebra $U_\eta(\L)$ has finite representation type and all its simple modules have dimension $M(\L)=p^{p^2-1}$.  
	\end{itemize}
However, the methods of \cite{PS99} failed when it came to showing that $U_\eta(\L)$ has fewer than $p^2$ pairwise non-isomorphic simple modules and hence is not semisimple (something that one would expect from a reduced enveloping algebra attached to a $T_0$-unstable $p$-character). It was clear that in dealing with such problems one would require a more detailed description of the centre of $U(\L)$. 
This unresolved issue is the starting point of the present investigation.

In Section~2, we describe the invariant algebra $S(\L)^\L$ which in characteristic $p>0$ contains
the large subalgebra $S(\L)^p=\{f^p\,|\,f\in S(\L)\}$. As $\L$ coincides with its derived subalgebra, $S(\L)^\L$ is a unique factorisation domain and hence is integrally closed in ${\rm Frac}\, S(\L)^\L=\k(\L^*)^\L$. We show that
the singular locus $\L^*_{\rm sing}:=\{\chi\in\L^*\,|\,\,\dim\L_\chi>2\}$ has codimension $\ge 2$ in $\L^*$ and produce two homogeneous functions
$F_1, F_2\in S(\L)^\L$ of degree $p+1$ and $p^2-p$ whose differentials 
are linearly independent at $\eta\in\L^*$. 

In order to establish the $\L$-invariance of $F_1$ we first show that it behaves nicely with respect to the action of $G$ on $S(\L)=\k[\L^*]$ (we use the standard convention that $g(f)=f\circ g^{-1}$ for all $f\in\k[\L^*]$ and $g\in G$). It turned out that in the notation of \ref{ss:1.1} we have
$$g(F_1)=\det(dg_0)^{p(p-2)}F_1 \ \ \ \ \mbox{and}\ \ \ g(F_2) =\det(dg_0)^{p(p-1)(p-2)}F_2\ \ \ \ (\forall\,g\in G).
$$
\noindent
Since the $G$-weight ${\rm wt}_G(F_1)$ is divisible by $p$ the Lie algebra $\L_{(0)}=\Lie(G)$ annihilates $F_1$. As $\L=\L_\eta\oplus \L_{(0)}$ and $\eta \not\in\L^*_{\rm sing}$ we now apply \cite[Corollary~5.3]{Sk} to deduce very quickly that $F_1\in S(\L)^\L$. For $F_2$ the argument is different: its $\L$-invariance is proved directly in Section~\ref{s:2}, and its semi-invariance is deduced in Section~\ref{s:3} from the first relation in \eqref{eq:1}. Note that  $$\deg(F_1)+\deg(F_2) = (p+1)+(p^2-p) =p^2+1= (\dim \L + {\rm ind}\, \L)/2,$$ 
the {\it magic number} of $\L$. Combined with the above-mentioned property of
$\L^*_{\rm sing}$ and the linear independence of  $(dF_1)_\eta$ and $(dF_2)_\eta$  this implies that
$$S(\L)^\L=S(\L)^p[F_1,F_2]$$ is a free $S(\L)^p$-module of rank $p^2$ with basis consisting of all
monomials $F_1^iF_2^j$ with $0\le i,j\le p-1$. 
\subsection{}
\label{ss:1.5} Let $U=U(\L)$ and write $U_d$ for the $d$-th component of the canonical filtration of $U$. If $d<p$ then the $(\ad\L)$-equivariant symmetrisation map $S^d(\L)\to  U_d$ allows one to lift symmetric invariants of degree $d$ to central elements of $U$ contained in $U_d$. 
This type of reasoning enabled N.N.~Yakovlev \cite{Yak} to determine, for $p>2$, the centre of $U(\g)$ in the case of the $p$-dimensional Witt algebra $\g=W(1)$. 
Unfortunately, such an approach would be unavailable in our case since $S(\L)^\L$ intersects trivially with $S^d(\L)$ for $0<d<p$.

Nevertheless, we show that all homogeneous elements of $S(\L)^\L$ can be lifted to central elements of $U(\L)$. The following is one of the main results of the paper:
\begin{Theorem}\label{main1}
	\label{m1}
There exist central elements $\tilde{F}_1\in Z\cap U_{p+1}$ and $\tilde{F}_2\in Z\cap U_{p^2-p}$ with principal symbols ${\rm gr}(\tilde{F_1})=F_1$ and ${\rm gr} (\tilde{F}_2)=F_2$. These elements  are semi-invariant under the action of $G$ on $U$, and  ${\rm wt}_G(\tilde{F}_i)={\rm wt}_G(F_i)$ for $i=1,2$.
\end{Theorem}
Using the results mentioned in \ref{ss:1.4} we deduce from Theorem~\ref{main1} that 
$$Z(U)\,=\,Z_p(\L)[\tilde{F}_1,\tilde{F}_2]\ \ \ \ \mbox{and}\ \ \ \ {\rm gr}\, Z(U)\,=\,S(\L)^\L.$$ Furthermore, $Z(U)$ is a free 
$Z_p(\L)$-module of rank $p^2$ with basis consisting of all monomials $\tilde{F}_1^i\tilde{F}_2^j$ with $0\le i,j\le p-1$.

We call the $\k$-subalgebra  generated by all central elements of $U$ semi-invariant under the action of $G$ the {\it Harish-Chandra centre} of $U$ and denote it by $Z_{\rm HC}(U)$. Arguing by induction on filtration degree and taking into account Theorem~\ref{main1} and the discussion in \ref{ss:1.4} one observes that
$$Z_{\rm HC}(U)\,=\,\k[\tilde{F}_1,\tilde{F}_2]$$ is a polynomial algebra in two variables.
\subsection{}
\label{ss:1.6}
Let $n=2p^2$. 
Recall from \ref{ss:1.3} that $Z_p(\L)$ is a polynomial algebra in  $z_i:=e_i^p-e_i^{[p]}$ where $\{e_1,\ldots, e_n\}$ is a $\k$-basis of $\L$. 
Given $f(\underline{e})=
f(e_1,\ldots, e_n)\in S(\L)$ we write $f(\underline{z})$ for the corresponding element of $Z_p(\L)$. To ease notation we denote the central elements $\tilde{F}_1$ and $\tilde{F}_2$ by $u$ and $v$, respectively. Our next result gives an explicit presentation of $Z$ by generators and relations:
\begin{Theorem}\label{main2}
\label{m2}
There is a semi-invariant $w\in S^{p+1}(\L)\setminus\{0\}$  of $G$-weight $p^2-p-1$
such that 
\begin{equation}u^p+uv=F_1(\underline{z})\ \ \ \mbox{and}\ \ \ v^p-w(\underline{z})^{p-2}u=F_2(\underline{z}).\label{eq:1}
\end{equation}
Moreover, ${\rm Specm}\,Z(U)$ is a complete intersection in the affine space
$\mathbb{A}^{n+2}(\k)$ with coordinates $u,v$ and $z_1,\ldots, z_n$ given by the equations \eqref{eq:1}.
\end{Theorem}
Recall the cocharacter $\mu\colon\,\k^\times\to T_0$ from \ref{ss:1.4}. 
Since $\det(d\mu(t)_0)=t^{p+1}$ and 
$$
F_i(\eta)=\mu(t)(F_i(\eta))=\big(\mu(t)(F_i)\big)(\mu(t)\cdot \eta)=t^{(p+1){\rm wt}_G(F_i)-	\deg(F_i)(p^2-p-1)}\cdot F_i(\eta)\ \ \ (\forall\, t\in\k^\times),$$ we see that $F_1(\eta)=F_2(\eta)=0$.
However,  such reasoning  would be inconclusive for $w\in S^{p+1}(\L)$ as it has $G$-weight $p^2-p-1$. In fact, we prove in  Section~\ref{s:3} that $w(\eta)=1$.
When specialised at $\eta$, the equations \eqref{eq:1} tell us that $u=v^p$ and $v^{p^2}+v^{p+1}=0$. 

Write $Z_\eta$ for the image of
$Z(U)$ in $U_\eta(\L)$ and set $S_\eta(\L):=S(\L)/\mathcal{I}_\eta$ where $\mathcal{I}_\eta$ is the ideal of $S(\L)$ generated by all 
$(x-\eta(x))^p$ with $x\in\L$. This algebra is the special member $U_{\eta,\,0}(\L)$ of the $\L$-equivariant
family  $\{U_{\eta,\, \lambda}(\L),  \,\, \lambda\in \k\}$  introduced in \cite[2.1]{PS99}.
From our discussion in \ref{ss:1.4} it follows that $U_{\eta,\,\lambda}(\L)\cong U_\eta(\L)$ for all nonzero $\lambda\in\k$. Let ${^\eta\!}{F}_i$ denote the image of $F_i$ in $S_\eta(\L)$.
Since $F_i(\eta)=0$ and $(dF_i)_\eta$ are linearly independent, the subalgebra of $S_\eta(\L)$ generated by  $^{\eta\!}F_1$ and  $^{\eta\!}F_2$ has dimension $p^2$. Furthermore, 
$S_\eta(\L)\,\cong\, \k[^{\eta\!}F_1, ^{\eta\!}F_2]\otimes S_\eta(\L_{(0)})$. Arguing in the spirit of \cite{PS99} and using Theorem~\ref{main1} we  deduce
that
$$U_\eta(\L)\,\cong\,Z_\eta\otimes U_\eta(\L_{(0)}).$$
This yields $\dim Z_\eta=p^2$ forcing $Z_\eta\cong \k[X]/\big(X^{p+1}(X^{p^2-p-1}+1)\big)$.

As $\L=\L_\eta\oplus \L_{(0)}$, the stabiliser of $\eta_{\,\vert \L_{(0)}}$ is trivial (hence toral). 
Applying \cite[Theorem~4.4(3)]{PS99} we get $U_\eta(\L_{(0)})\,\cong {\rm Mat}_N(\k)$ where $N=p^{p^2-1}$. The Chinese Remainder Theorem entails that $$U_\eta(\L)\,\cong\,{\rm Mat}_N(Z_\eta) \,\cong\,{\rm Mat}_N\big(\k[X]/(X^{p+1})\big)\oplus \underbrace{{\rm Mat}_N(\k)\oplus\cdots\oplus {\rm Mat}_N\big(\k)}_{p^2-p-1 \text{ copies}}$$
(one should keep in mind that $p^2-p-1$ is coprime to the characteristic of $\k$).
We thus obtain that $U_\eta(\L)$ has $p^2-p-1$ simple blocks and one local block Morita equivalent 
to $\k[X]/(X^{p+1})$. As it turned out, the simple blocks form a single orbit under the action of the cyclic group
$\langle \mu(\zeta)\rangle \subset {\rm Aut}\,U_\eta(\L)$ where $\zeta$ is a primitive  $(p^2-p-1)$-st root of unity in $\k$. In all, the reduced enveloping algebra $U_\eta(\L)$ has $p^2-p$ pairwise non-isomorphic simple modules, which solves the problem left open in \cite[6.4]{PS99}.
\subsection{}
\label{ss:1.7} In Section~\ref{s:4}, we prove that the set $$\mathcal{X}_{\rm sing}\,:=\,\{\chi\in\L^*\,\vert\,\, U_\chi(\L)\  \mbox{has a simple module of dimension}\,<N\}$$ is contained in $\L^*_{\rm sing}$, which implies that
$\Psi^{-1}(\L^*_{\rm reg})\subseteq \Omega$. Our proof 
relies on a careful study of the open set $\L^*_{\rm reg}:=\L^*\setminus\L^*_{\rm sing}$ and 
makes use of \cite[Proposition~5.2]{PS99}. 
The explicit presentation of $Z(U)$ provided by Theorem~\ref{m2} enables us to determine the singular points of its maximal spectrum. We then use the fact that the Azumaya locus 
$\Omega\subseteq {\rm Specm}\,Z(U)$ coincides with the smooth locus of $Z(U)$ to show that
$\Psi^{-1}(\chi)\cap\Omega=\varnothing$ for all $\chi\in \L^*_{\rm sing}$. As a corollary, we get $\mathcal{X}_{\rm sing}=\L^*_{\rm sing}$ and obtain our third main result:
\begin{Theorem}\label{main3}
The Azumaya locus $\Omega$ of $Z(U)$ coincides with the inverse image of $\L^*_{\rm reg}$
under the finite morphism $\Psi\colon\, {\rm Specm}\,Z(U)\twoheadrightarrow \L^*$.	
\end{Theorem}
\noindent
In particular, $\Psi(\Omega)=\L^*_{\rm reg}$. Corollary~\ref{C:4.6} gives an explicit description of this set in terms of normal forms of residue duals (to be defined in Section~\ref{s:2}), which settles the case $n=2$ of the problem mentioned in \ref{ss:1.3}.

\noindent
At the end of the paper we show that the geometric monodromy group of the finite morphism
${\rm Specm}\,Z(U)\twoheadrightarrow {\rm Specm}\,Z_p(\L)$ is isomorphic to the symmetric group $\mathfrak{S}_{p^2}$ and then use the Galois correspondence  to show  that the field extension $Q(\L)/Q_p(\L)$ has no intermediate subfields $E/Q_p(\L)$. We also explain why Theorem~\ref{main3} is no longer valid in characteristics $2$ and $3$.

\smallskip
\noindent
{\bf Acknowledgment.} After this preprint first appeared on arXiv, we were informed by Leonid Bedratyuk that he had independently proved that $S(\L)^\L = S(\L)^p[F_1, F_2]$ and
$Z(U) = Z_p(\L)[u, v]$. This appears in his preprint \cite{Bed26}, which relies on methods different from ours. His preprint covers the characteristic $3$ case as well, but it does not give a presentation of $Z(U)$ (our Theorem~\ref{main2}) and does not discuss the Azumaya locus of $Z(U)$.
We note that symmetric invariants of $W(n)$ were studied earlier in \cite{Bed93, Bed94}.
\section{Symmetric invariants of $\L$}
\label{s:2}
\subsection{}
\label{ss:2.1}
Let $\Lambda_p:=\{(a,b)\in \Z^2\,\vert\,\,0\le a,b\le p-1\}.$
Given $f\in\O$ we put ${\rm Res}(f):=\del_x^{p-1}\del_y^{p-1}(f)$, the coefficient of $f$ at $x^{p-1}y^{p-1}$.
 One of our tools in this paper will be the symmetric bilinear form  on $\O$  which assigns to $(f,g)\in \O\times \O$ the scalar ${\rm Res}(fg)$.
This form is non-degenerate as
for $(a,b), (c,d)\in\Lambda_p$ we have  ${\rm Res}(x^a y^b\cdot  x^c y^d)=0$ if $(c,d)\neq (p-1-a,p-1-b)$ and ${\rm Res}(x^ay^b \cdot x^{p-1-a}y^{p-1-b})=1$. 

Let $X=X(x)\del_x+X(y)\del_y$ and $Y=Y(x)\del_x+Y(y)\del_y$ be two arbitrary elements of $\L$ and set
\begin{equation}
	\nu(X,Y)\,:=\,X(x)Y(y)-X(y)Y(x)\in \O.\label{eq:2.1}
\end{equation}
Composing ${\rm Res} \in \O^*$ with $\nu$ we obtain a non-degenerate alternating bilinear form 
$$\omega(X,Y)\,:=\,{\rm Res}\big(\nu(X,Y)\big)\ \ \ \ \ \ \ \ \ \ \ \ (\forall\,X,Y\in \L)$$
on $\L$ which enables us to identify $\L^*$ and $\L$ as vector spaces: for every $\xi\in\L^*$ there is a unique $J\xi\in \L$ such that $\xi(Y)=\omega(J\xi,Y)$ for all $Y\in \L$. We call $J\xi$ the {\it residue dual} of $\xi$.\footnote{The pairing $\omega$ is the residue of the $2$-form $X\wedge Y$ contracted with $dx\wedge dy$; this is the classical identification of $W(2)^*$ with $2$-densities, cf.\ the module $\L_{2\,{\rm div}}$ below.} 

For $D=f\del_x+g\del_y\in \L$ we denote by ${\rm div}\,D$ the {\it divergence} of $D$. It is given by the formula ${\rm div}(D)=\del_x(f)+\del_y(g)$ and has the property that 
$${\rm div}([X,Y]) = X({\rm div}(Y))-Y({\rm div}(X))$$
for all $X,Y\in \L$; see \cite[2.9]{BGP}, for example. This enables one to define a new module structure on $\L$ by setting
$$D\cdot X:=[D,X]+2{\rm div}(D)X\ \ \ \ \ \ \ (\forall\,D,X\in \L).$$  
The $\L$-module thus obtained will be denoted by $\L_{2\,{\rm div}}$. Since $\L$ is a simple restricted $(\ad \L)$-module of dimension $2p^2$, so is the coadjoint $\L$-module $\L^*$. In conjunction with our discussion at the end of \ref{ss:1.1} this implies that $\L^*$ is isomorphic to the induced module $U_0(\L)\otimes_{U_0(\L_{(0)})}\,\L^*_{-1}$, and $\L^*_{-1}$ is an irreducible $\L_0$-module of highest weight $(1,0)$ with respect to $\Lie(B_0)$.

Next we observe that the $\L$-module $\L_{2\,\rm div}$  has a submodule generated by the subspace $\L_{2p-3}$ which is still invariant under the action of $\L_{(0)}$. As ${\rm div}(x\del_x)={\rm div}(y\del_y)=1$ the vector $x^{p-1}y^{p-1}\del_y$ of $\L_{2\,{\rm div}}$ 
acquires highest weight $ (-1,-2)+(2,2)=(1,0)$ with respect to $\Lie(B_0)$. Thanks to the irreducibility of the induced module $\L^*$ this entails that  $\L_{2\,{\rm div}}\cong \L^*$ as $\L$-modules. We stress that this isomorphism is specific to $W(2)$ and would require modifications to work for $W(n)$ with $n\ge 3$.
\begin{Lemma}[\rm Residue twist] \label{L:2.1}
	The map $J\colon\,\L^*\to \L_{2\,{\rm div}}$ is an isomorphism of $\L$-modules, that is, $J(D\cdot \xi)=[D,J\xi]+2{\rm div}(D)J\xi$ for all $\xi\in\L^*$ and $D\in\L$.
	\end{Lemma}
	\begin{proof}
For $X\in\L$ let $\vec X=\big(X(x),X(y)\big)^{\sf T}\in\O^2$ be its column of coefficients, so that $\nu(X,Y)=\det\big(\vec X\ \vec Y\big)$ is the determinant of the $2\times2$ matrix with columns $\vec X$ and $\vec Y$. For $D\in\L$ write $D(\vec X)=\big(D(X(x)),D(X(y))\big)^{\sf T}$, and let
$$A_D=\begin{pmatrix}\del_x(D(x))&\del_y(D(x))\\ \del_x(D(y))&\del_y(D(y))\end{pmatrix}\in{\rm Mat}_2(\O),$$
so that ${\rm tr}(A_D)={\rm div}(D)$. Since ${\rm Res}(\del_x(\O))={\rm Res}(\del_y(\O))=0$, the product rule gives ${\rm Res}(D(a))=-{\rm Res}\big({\rm div}(D)a\big)$ for all $a\in\O$.

As $D$ is a derivation of $\O$, applying it to the determinant gives
$$D\big(\nu(X,Y)\big)=\det\big(D(\vec X)\ \vec Y\big)+\det\big(\vec X\ D(\vec Y)\big).$$
For $c\in\{x,y\}$ we have $[D,X](c)=D(X(c))-X(D(c))$, implying that
$$\overrightarrow{[D,X]}=D(\vec X)-A_D\vec X.$$
Every $A\in{\rm Mat}_2(\O)$ satisfies $\det(A\vec X\ \vec Y)+\det(\vec X\ A\vec Y)={\rm tr}(A)\det(\vec X\ \vec Y)$. Substituting $D(\vec X)=\overrightarrow{[D,X]}+A_D\vec X$ and $D(\vec Y)=\overrightarrow{[D,Y]}+A_D\vec Y$ we obtain
$$D\big(\nu(X,Y)\big)=\nu([D,X],Y)+\nu(X,[D,Y])+{\rm div}(D)\,\nu(X,Y).$$
Applying ${\rm Res}$ and using ${\rm Res}(D(a))=-{\rm Res}\big({\rm div}(D)a\big)$ with $a=\nu(X,Y)$, we get
$$\omega([D,X],Y)+\omega(X,[D,Y])=-2\,\omega\big({\rm div}(D)X,Y\big)\qquad(\forall\,X,Y\in\L).$$
Now let $X=J\xi$. Then
$$(D\cdot\xi)(Y)=-\xi([D,Y])=-\omega(X,[D,Y])=\omega\big([D,X]+2{\rm div}(D)X,\,Y\big)$$
for all $Y\in\L$. As $\omega$ is non-degenerate, $J(D\cdot\xi)=[D,J\xi]+2{\rm div}(D)J\xi$.
\end{proof}
\subsection{}
\label{ss:2.2}
For $f\in S(\L)=\k[\L^*]$ and $\xi\in\L^*$ we denote by $(df)_\xi\in(\L^*)^*=\L$ the differential of $f$ at $\xi$, so that $f(\xi+\epsilon\zeta)=f(\xi)+\epsilon\,\zeta\big((df)_\xi\big)$ in $\k[\epsilon]/(\epsilon^2)$ for all $\zeta\in\L^*$. The algebra $\L$ acts on $S(\L)$ by derivations, and $(D\cdot f)(\xi)=-(D\cdot\xi)\big((df)_\xi\big)$; thus $f\in S(\L)^\L$ if and only if $(D\cdot\xi)\big((df)_\xi\big)=0$ for all $D\in\L$ and $\xi\in\L^*$.

For $X\in\L$ we denote by $\psi(X)$ the coefficient of $T^p$ in the characteristic polynomial $\det(T\,{\rm Id}-X)$ of the linear map $X\colon\O\to\O$. It is proved in \cite{Pr92} that
$\det(T\,{\rm Id}-X)$ is a $p$-polynomial in $T$ of degree $p^2$ whose non-constant coefficients, viewed as functions on $\L$, freely generate the invariant ring $\k[\L]^G$.
Recall that $X^{[p]}=X^p$, and define polynomial functions $F_1$ and $F_\psi$ on $\L^*$ by setting
\begin{equation}\label{eq:2.2}
	F_1(\xi)\,:=\,\omega\big(J\xi,(J\xi)^{[p]}\big)\ \ \ \ \mbox{and}\ \ \ \ F_\psi(\xi)\,:=\,\psi(J\xi)\ \ \ \ \ \ (\forall\,\xi\in\L^*).
\end{equation}
They are homogeneous of degree $p+1$ and $p^2-p$, respectively, and have coefficients in $\mathbb{F}_p$ with respect to the monomial basis of $\L$.

We need two properties of the characteristic polynomial. For $X,H\in\End_\k(\O)$, the derivative of $X\mapsto\det(T\,{\rm Id}-X)$ at $X$ in the direction $H$ is
\begin{equation}\label{eq:2.3}
	d\det(T\,{\rm Id}-X)(H)\,=\,-{\rm tr}\big({\rm adj}(T\,{\rm Id}-X)\,H\big),
\end{equation}
where ${\rm adj}$ denotes the adjugate matrix. Moreover,
\begin{equation}\label{eq:2.4}
	(d\psi)_X(aX)\,=\,0\ \ \ \ \ \ (\forall\,X\in\L,\ a\in\O).
\end{equation}
Indeed, the left hand side of \eqref{eq:2.4} is polynomial in $X$ and $a$, so we may assume that $X=\lambda_1(1+x)\del_x+\lambda_2(1+y)\del_y$ with $\lambda_1,\lambda_2$ linearly independent over $\mathbb{F}_p$. Then $X$ has pairwise distinct eigenvalues on the basis $\{(1+x)^r(1+y)^s\,|\,\,(r,s)\in\Lambda_p\}$ of $\O$, its zero eigenspace is $\k1$, and the product of two eigenvectors is an eigenvector whose eigenvalue is the sum of theirs. Hence the diagonal part of multiplication by $a$ in this basis is $a(0){\rm Id}$, and \eqref{eq:2.3} together with Euler's identity gives $(d\psi)_X(aX)=a(0)(p^2-p)\psi(X)=0$.
\begin{Proposition}\label{P:2.2}
The functions $F_1$ and $F_\psi$ lie in $S(\L)^\L$.
\end{Proposition}
\begin{proof}
Let $D\in\L$, $\xi\in\L^*$, and put $X=J\xi$, $a={\rm div}(D)$ and $\zeta=D\cdot\xi$. By Lemma~\ref{L:2.1}, $J\zeta=[D,X]+2aX$. For $f=F_1$ we have $\zeta\big((df)_\xi\big)=\frac{d}{d\epsilon}F_1(\xi+\epsilon\zeta)\big|_{\epsilon=0}=\omega(J\zeta,X^p)+\omega\big(X,(\ad X)^{p-1}(J\zeta)\big)$, because the derivative of $X\mapsto X^p$ in the direction $H$ is $\sum_{i=0}^{p-1}X^iHX^{p-1-i}=(\ad X)^{p-1}(H)$. Since $(\ad X)^{p-1}([D,X])=[D,X^p]$, the contribution of $[D,X]$ is $\omega([D,X],X^p)+\omega(X,[D,X^p])=-2\omega(aX,X^p)$ by the identity established in the proof of Lemma~\ref{L:2.1}. The contribution of $2aX$ is $2\omega(aX,X^p)+2\omega\big(X,(\ad X)^{p-1}(aX)\big)$. Since $(\ad X)^r(aX)=X^r(a)X$ for all $r\ge 1$ and $\omega$ is alternating, the second term vanishes. Adding up we obtain $\zeta\big((dF_1)_\xi\big)=0$, so $F_1\in S(\L)^\L$.

For $f=F_\psi$ we have $\zeta\big((dF_\psi)_\xi\big)=(d\psi)_X(J\zeta)$. Since ${\rm adj}(T\,{\rm Id}-X)$ is a polynomial in $X$ and the trace of $P(X)[D,X]$ is zero for every polynomial $P$, formula \eqref{eq:2.3} gives $(d\psi)_X([D,X])=0$; moreover, $(d\psi)_X(2aX)=0$ by \eqref{eq:2.4}. Therefore $F_\psi\in S(\L)^\L$.
\end{proof}

\subsection{}
\label{ss:2.3}
Recall that $\Delta=\del_x-x^{p-1}\del_y$ and $\eta(Y)=\omega(\Delta,Y)$ for all $Y\in\L$. The monomials $x^ay^b$ with $(a,b)\in\Lambda_p$, ordered by the value of $a+pb$, form a single nilpotent chain for $\Delta$: one has $\Delta(x^ay^b)=ax^{a-1}y^b$ if $a>0$ and $\Delta(y^b)=-by^{b-1}x^{p-1}$. Using Wilson's theorem one obtains
\begin{equation}\label{eq:2.5}
	\Delta^p=\del_y,\ \ \ \ \ \ \Ker\Delta=\k1,\ \ \ \ \ \ \Delta^{p^2-1}(f)={\rm Res}(f)\ \ \ \ \ (\forall\,f\in\O).
\end{equation}
Similar results hold for all Witt--Jacobson Lie algebras $W(n)$; see \cite{Pr92, Pr15} for detail.
\begin{Lemma}\label{L:2.3}
One has $\L_\eta=\k\Delta\oplus\k\Delta^{[p]}$, and the differentials of $F_1$ and $F_\psi$ at $\eta$ are $(dF_1)_\eta=\del_y$ and $(dF_\psi)_\eta=\Delta$. In particular, ${\rm ind}\,\L=2$.
\end{Lemma}
\begin{proof}
Let $Y=a\Delta+b\del_y$ with $a,b\in\O$. Since ${\rm div}(\Delta)=0$ and $\Delta^{[p]}=\del_y$ commutes with $\Delta$, Lemma~\ref{L:2.1} shows that $Y\in\L_\eta$ if and only if
$$\big(\Delta(a)+2\del_y(b)\big)\Delta-\Delta(b)\,\del_y\,=\,0.$$
As $\Delta$ and $\del_y$ form an $\O$-basis of $\L$, it follows from \eqref{eq:2.5} that $b\in\k$, and then $a\in\k$. Since ${\rm div}(\Delta)=0$, differentiating \eqref{eq:2.2} at $\eta$ gives $(dF_1)_\eta(\zeta)=\zeta(\del_y)$ for all $\zeta\in\L^*$. For $Y=f\Delta+g\del_y$, formula \eqref{eq:2.3} yields
$$(d\psi)_\Delta(Y)\,=\,-{\rm tr}\big(\Delta^{p^2-p-1}Y\big)\,=\,-{\rm Res}(g),$$
because ${\rm tr}(\del_y^{p-1}f)=0$ (the $p$ equal diagonal contributions add up to zero) and ${\rm tr}(\Delta^{p^2-1}g)={\rm Res}(g)$ by \eqref{eq:2.5}. Hence $(dF_\psi)_\eta=\Delta$. For $f\in S(\L)^\L$ and $D\in\L$ we have $\xi([D,(df)_\xi])=-(D\cdot\xi)((df)_\xi)=0$, so $(df)_\xi\in\L_\xi$ for all $\xi$. Since $(dF_1)_\xi$ and $(dF_\psi)_\xi$ are linearly independent at $\xi=\eta$, they are linearly independent on a nonempty open subset of $\L^*$, where therefore $\dim\L_\xi\ge2$. By upper semicontinuity of $\xi\mapsto\dim\L_\xi$, we get ${\rm ind}\,\L\ge2$, and $\dim\L_\eta=2$ gives equality.
\end{proof}

\subsection{}
\label{ss:2.4}
For $X\in\L$ put $q_X:=\big(X+2{\rm div}(X)\big)^{p-1}(1)\in\O$, where ${\rm div}(X)$ acts on $\O$ by multiplication, and define
$$h_0(\xi)\,:=\,q_{J\xi}(0),\qquad w(\xi)\,:=\,\nu\big(J\xi,(J\xi)^{[p]}\big)(0)\ \ \ \ \ \ (\forall\,\xi\in\L^*).$$
These are homogeneous polynomial functions on $\L^*$ of degree $p-1$ and $p+1$, respectively, with coefficients in $\mathbb{F}_p$. We set
$$F_2\,:=\,F_\psi+h_0^p.$$
Since $h_0^p\in S(\L)^p$, Proposition~\ref{P:2.2} shows that $F_2\in S(\L)^\L$, and $(dF_2)_\eta=(dF_\psi)_\eta=\Delta$.

For $g\in G$ and $X\in\L$ put $g_*(X)=g\circ X\circ g^{-1}\in\L$, and let $j_g=\det\big(\del_j(g(x_i))\big)_{1\le i,j\le2}\in\O^\times$, where $(x_1,x_2)=(x,y)$, so that $j_g(0)=\det(dg_0)$. One has ${\rm Res}\big(g(f)j_g\big)=\det(dg_0)^p\,{\rm Res}(f)$ for all $f\in\O$.\footnote{The linear form $f\mapsto{\rm Res}(f)$ on $\O\,dx\wedge dy$ vanishes exactly on the exact forms $d(\O\,dx)+d(\O\,dy)$, and $g$ maps $f\,dx\wedge dy$ to $g(f)j_g\,dx\wedge dy$; hence $f\mapsto{\rm Res}(g(f)j_g)$ is a multiple of ${\rm Res}$, and evaluating on $f=x^{p-1}y^{p-1}$ gives the factor $\det(dg_0)^{p}$.} This implies that
\begin{equation}\label{eq:2.6}
	J(g\cdot\xi)\,=\,\det(dg_0)^{-p}j_g^2\,g_*(J\xi)\ \ \ \ \ \ (\forall\,\xi\in\L^*).
\end{equation}
\begin{Lemma}\label{L:2.4}
For all $g\in G$ we have $g(F_1)=\det(dg_0)^{p(p-2)}F_1$ and $g(w)=\det(dg_0)^{p^2-p-1}w$.
\end{Lemma}
\begin{proof}
Put $\delta=\det(dg_0)$, $X=J\xi$ and $f=\delta^{-p}j_g^2$, so that $f(0)=\delta^{2-p}$. For $Y\in\L$ and $f\in\O$ one has $(fY)^{[p]}=f(0)^pY^{[p]}+(fY)^{p-1}(f)\,Y$ (Hochschild's formula, proved in the first paragraph of the proof of Lemma~\ref{L:3.2}; note $f^p=f(0)^p$ in $\O$), and hence $\nu\big(fY,(fY)^{[p]}\big)=f(0)^pf\,\nu(Y,Y^{[p]})$. On the other hand, $\nu(g_*Y,g_*Y')=j_g^{-1}g\big(\nu(Y,Y')\big)$ for all $Y,Y'\in\L$. Combining this with \eqref{eq:2.6} we get
$$\nu\big(J(g\cdot\xi),J(g\cdot\xi)^{[p]}\big)\,=\,\delta^{p(2-p)-p}\,j_g\,g\big(\nu(X,X^{[p]})\big).$$
Evaluating at the origin gives $w(g\cdot\xi)=\delta^{-(p^2-p-1)}w(\xi)$, whereas applying ${\rm Res}$ gives $F_1(g\cdot\xi)=\delta^{-p(p-2)}F_1(\xi)$. Since $g(f)=f\circ g^{-1}$, the lemma follows.
\end{proof}
The function $F_2$ is also a $G$-semi-invariant, of weight $p(p-1)(p-2)$; this will be deduced in Lemma~\ref{L:3.7} from the first relation in \eqref{eq:1}.

\subsection{}
\label{ss:2.5a}
We now describe the linear functions $\chi$ with $J\chi\notin\L_{(0)}$ up to the action of $G$; this yields the codimension estimate for $\L^*_{\rm sing}$ and will be used again in Section~\ref{s:4}.
\begin{Lemma}\label{L:4.2}
If $J\chi\notin\L_{(0)}$ then, after interchanging $x$ and $y$ if necessary, $\chi$ is $G$-conjugate to a linear function whose residue dual has the form
\begin{equation}\label{eq:4.2}
	f\cdot\big(\del_x+x^{p-1}c(y)\del_y\big),\qquad f\in\O^\times,\ \ c\in\k[y]/(y^p).
\end{equation}
Such a linear function is regular if and only if $c(0)\ne 0$ or $c'(0)\ne 0$.
\end{Lemma}
\begin{proof}
Suppose $X=J\chi$ has $X(x)\in\O^\times$ and set $e=X(x)^{-1}X$. Then $e=\del_x+\sum_{j=0}^{p-1}x^j\Gamma_j$ with $\Gamma_j\in\Der\,\k[y]/(y^p)$. If the coefficients $\Gamma_0,\ldots,\Gamma_{j-1}$ have already been removed and $j\le p-2$, put $H_j=-\frac{x^{j+1}}{j+1}\Gamma_j$ and $\alpha_j=\exp(-H_j)$. The derivation $H_j$ fixes $x$ and raises the $x$-degree by at least $j+1$; hence $H_j^p=0$ and $\alpha_j$ is an automorphism of $\O$. Modulo $x^{j+1}\Der\O$ we have $(\alpha_j)_*e\equiv e-[H_j,e]\equiv\del_x$, because $[\del_x,H_j]=-x^j\Gamma_j$ and all other terms have $x$-degree $>j$. After $p-1$ steps we arrive at $\del_x+x^{p-1}c(y)\del_y$. Each step changes the later coefficients by polynomials in the earlier ones, so the eliminated coefficients together with $c$ form affine coordinates. By \eqref{eq:2.6}, the coadjoint action differs from conjugation by a unit factor, which gives \eqref{eq:4.2}.

Let $X=f\delta$ with $\delta=\del_x+x^{p-1}C$, $C=c(y)\del_y$, and write $D=a\delta+E$ with $E(x)=0$. By Lemma~\ref{L:2.1}, $D\in\L_\chi$ if and only if $[D,f\delta]+2{\rm div}(D)f\delta=0$. Expanding $E=\sum_jx^jE_j$ with $E_j\in\Der\,\k[y]/(y^p)$, the transverse component gives $E_j=0$ for $j\ge 1$ and $[E_0,C]=0$. Writing $b=fa$, the remaining component reads
\begin{equation}\label{eq:4.3}
	\big(\del_x+x^{p-1}T_c\big)(b)\,=\,-E(f)-2E'f,\qquad T_c(g)=cg'+2c'g,
\end{equation}
where $E=E_0=e(y)\del_y$ and $E'=e'(y)$. For fixed $E$ the solutions form an empty set or a coset of $\Ker T_c$, so that $\dim\L_\chi\le\dim\Ker T_c+\dim\mathfrak{c}(C)$, where $\mathfrak{c}(C)$ is the centraliser of $C$ in $\Der\,\k[y]/(y^p)$. If $c$ is a unit, then $(c^2g)'=cT_c(g)$ and $ec'-ce'=-c^2(e/c)'$ show that both spaces are one-dimensional. If $c=c_1y+\cdots$ with $c_1\ne 0$, their matrices in the basis $1,y,\ldots,y^{p-1}$ are triangular with diagonal entries $c_1(j+2)$ and $c_1(1-j)$, so both kernels are again one-dimensional. As ${\rm ind}\,\L=2$, the linear function is regular in both cases.

Conversely, $c=0$ gives $p$ independent solutions of \eqref{eq:4.3} with $E=0$. If $s={\rm ord}_y(c)\ge 3$, the monomials $y^{p-s+1},\ldots,y^{p-1}$ lie in $\Ker T_c$, and if $s=2$ the matrix of $T_c$ from $1,\ldots,y^{p-2}$ to $y,\ldots,y^{p-1}$ is triangular with diagonal entries $c_2(j+4)$, one of which vanishes; in addition $y^{p-1}\in\Ker T_c$. Hence $\dim\Ker T_c\ge 2$ whenever $c\ne 0$ and $s\ge 2$. For such $c$, writing $f=\sum_ix^if_i(y)$, the element $b=-f_{p-1}-\sum_{i=1}^{p-1}\frac{x^i}{i}T_c(f_{i-1})$ solves \eqref{eq:4.3} with $E=C$. This solution is independent of those with $E=0$, so that $\dim\L_\chi\ge 3$, and hence $\dim\L_\chi\ge 4$, because $\L_\chi$ is the radical of the alternating form $(X,Y)\mapsto\chi([X,Y])$ on the even-dimensional space $\L$.
\end{proof}
Lemma~\ref{L:4.2} shows that $\L^*_{\rm sing}$ has codimension $\ge 2$: on each of the two charts $J\chi(x)\in\O^\times$ and $J\chi(y)\in\O^\times$ the coefficients $c(0)$ and $c'(0)$ are independent affine coordinates, while the complement of these charts has codimension two.

\subsection{}
\label{ss:2.5}
The set $\L^*_{\rm sing}$ has codimension $\ge 2$ in $\L^*$; see \ref{ss:2.5a} and \cite[6.4]{PS99}. Choose a basis $e_1,\ldots,e_n$ of $\L$, where $n=2p^2$, and let $\Pi_{ij}=(-1)^{i+j}{\rm Pf}\big(M_{\widehat{ij}}\big)$ for $i<j$, where $M=\big(\xi([e_i,e_j])\big)_{1\le i,j\le n}$ and $M_{\widehat{ij}}$ is obtained from $M$ by deleting the rows and columns with indices $i$ and $j$. Put $\Pi=\sum_{i<j}\Pi_{ij}\,e_i\wedge e_j$.
\begin{Lemma}\label{L:2.5}
There is $\gamma\in\k^\times$ such that $dF_1\wedge dF_2=\gamma\,\Pi$.
\end{Lemma}
\begin{proof}
The Pfaffian expansion shows that the matrix $(\Pi_{ij})$ annihilates $\big(\xi([e_i,e_j])\big)$. On $\L^*_{\rm reg}$ the matrix $(\Pi_{ij})$ has rank two and its image is $\L_\xi$, so that $\Pi$ spans $\bigwedge^2\L_\xi$. All $\Pi_{ij}$ vanish on $\L^*_{\rm sing}$, and since this set has codimension $\ge 2$, they have no nonconstant common divisor. The differentials of $F_1$ and $F_2$ lie in $\L_\xi$ and are linearly independent at $\eta$ by Lemma~\ref{L:2.3}. Therefore $dF_1\wedge dF_2=q\Pi$ for some $q\in\k(\L^*)$. The denominator of $q$ divides all $\Pi_{ij}$, hence $q$ is a polynomial. Both sides have degree $p^2-1$, forcing $q\in\k^\times$.
\end{proof}
Our next result is a direct consequence of \cite[Theorem~5.4]{Sk}. We include a short proof for the reader's convenience.
\begin{Theorem}\label{T:2.6}
$S(\L)^\L=S(\L)^p[F_1,F_2]$ is a normal complete intersection and a free $S(\L)^p$-module with basis $\{F_1^iF_2^j\,|\,\,0\le i,j\le p-1\}$.
\end{Theorem}
\begin{proof}
Put $R=S(\L)$ and $K={\rm Frac}\,R$. The coadjoint vector fields span an $(n-2)$-dimensional subspace of $\Der_{K^p}K$ at the generic point, and the differentials of a $p$-independent family in $K^\L$ (that is, a family whose differentials in $\Omega_{K/K^p}$ are linearly independent) annihilate this subspace. Hence such a family has at most two members and $[K^\L\colon K^p]\le p^2$. As $(dF_1)_\eta$ and $(dF_2)_\eta$ are linearly independent, $F_1$ and $F_2$ are $p$-independent, and we get $K^\L=K^p(F_1,F_2)$.

Let $B=R^p[F_1,F_2]$. Since $F_1,F_2$ are $p$-independent and $F_i^p\in R^p$, the monomials $F_1^iF_2^j$ with $0\le i,j\le p-1$ form a basis of $B$ over $R^p$. Writing $R=\k[t_1,\ldots,t_n]$ and $y_i=t_i^p$, we obtain
$$B\,\cong\,\k[y_1,\ldots,y_n,s_1,s_2]/\big(s_1^p-\overline{F}_1(y),\,s_2^p-\overline{F}_2(y)\big),$$
where $\overline{F}_i(t_1^p,\ldots,t_n^p)=F_i(t_1,\ldots,t_n)^p$. This is a complete intersection, hence satisfies Serre's condition $(S_2)$. The rows of its Jacobian matrix are $-d\overline{F}_1$ and $-d\overline{F}_2$, and $\big(\del\overline{F}_i/\del y_j\big)(t^p)=\big(\del F_i/\del t_j\big)(t)^p$. By Lemma~\ref{L:2.5}, the Jacobian matrix has rank $<2$ exactly over $\L^*_{\rm sing}$. The morphism ${\rm Spec}\,B\to{\rm Spec}\,R^p$ is finite and bijective on points, since $B^p\subseteq R^p$, and $R^p\cong R$ via the Frobenius map. Hence it identifies closed subsets and preserves codimension. Since the singular locus of ${\rm Spec}\,B$ has codimension $\ge 2$, Serre's criterion shows that $B$ is normal. Every element of $R^\L$ is integral over $B$ and lies in $K^\L={\rm Frac}\,B$, hence belongs to $B$.
\end{proof}
Since $F_1$ and $F_2$ are annihilated by $\L$, the weights in Lemma~\ref{L:2.4} and Lemma~\ref{L:3.7} are divisible by $p$, in agreement with the fact that $\L_{(0)}=\Lie(G)$ acts trivially on $S(\L)^\L$.

\section{The centre of $U(\L)$}
\label{s:3}
\subsection{}
\label{ss:3.1}
Throughout this section we use the notation of Section~\ref{s:2}. Thus $\eta(Y)=\omega(\Delta,Y)$ for all $Y\in\L$, the functions $F_1,F_\psi,h_0,w$ and $F_2=F_\psi+h_0^p$ are defined in \ref{ss:2.2} and \ref{ss:2.4}, and $S(\L)^\L=S(\L)^p[F_1,F_2]$ by Theorem~\ref{T:2.6}. Since $h_0^p\in S(\L)^p$, the latter statement holds equally with $F_\psi$ in place of $F_2$.

We fix the monomial basis $\{e_b\}$ of $\L$ consisting of the vector fields $x^iy^j\del_x$ and $x^iy^j\del_y$, and denote by $\{g_b\}$ the basis of $\L$ such that $\omega(g_b,e_c)=\delta_{bc}$. As in \ref{ss:1.6} we write $z_b=e_b^p-e_b^{[p]}$ and, for a polynomial $f$ with coefficients in $\mathbb{F}_p$, we denote by $f(\underline{z})$ the corresponding element of $Z_p(\L)$. For $a\in U_d\setminus U_{d-1}$ its {\it principal symbol} is the image of $a$ in $U_d/U_{d-1}=S^d(\L)$; we write it ${\rm gr}(a)$. The principal symbol of $f(\underline{z})$ is $f^p$. The standard grading of $\L$ induces a $\Z$-grading of $U$ and $S(\L)$ which we call the {\it integer grading}. If $e_b$ has degree $d_b$ then $g_b$ has degree $2p-4-d_b$, and hence $F_1$ and $F_\psi$ are homogeneous of integer degrees
\begin{equation}\label{eq:3.1}
	\theta_1=2p(p-2)\ \ \ \ \mbox{and}\ \ \ \ \theta_2=(p-1)\theta_1.
\end{equation}

\subsection{}
\label{ss:3.2}
We first lift $F_1$. Consider the algebra $\End_\k(\O)\otimes U$ and its element\footnote{$\D$ is the image of the canonical element of $\L\otimes\L$ under the identification of $\L^*$ with $\L$ given by $\omega$. Replacing $e_b$ by $\xi(e_b)$ turns $\D$ into the vector field $J\xi$; the element $u$ below is a noncommutative version of $\xi\mapsto\omega(J\xi,(J\xi)^{[p]})$.}
$$\D\,:=\,\textstyle{\sum_b}\,g_b\otimes e_b\,=\,a_1\del_x+a_2\del_y,\ \ \ \ \ \ a_1=\textstyle{\sum_b}\,g_b(x)\otimes e_b,\ \ \ a_2=\textstyle{\sum_b}\,g_b(y)\otimes e_b,$$
where all coefficients act on $\O\otimes U$ by left multiplication. For $\Phi\in \End_\k(\O)\otimes U$ we set $\R(\Phi):={\rm Res}\big(a_1\Phi(y)-a_2\Phi(x)\big)\in U$ and define
$$u\,:=\,\R(\D^p).$$
\begin{Lemma}\label{L:3.1}
{\rm (i)} If $|\alpha|+|\beta|\le 2p-4$ then $[\del^\alpha a_i,\del^\beta a_j]=0$ for $i,j\in\{1,2\}$.

\smallskip

\noindent
{\rm (ii)} The map $\O\to \O\otimes U$, $f\mapsto \D^p(f)$, is a derivation.
\end{Lemma}
\begin{proof}
Write $X_1,X_2$ for the coordinates in the copy of $\O$ on which $\End_\k(\O)$ acts and $x,y$ for those in $\L$. The dual vectors of $x^cy^d\del_x$ and $x^cy^d\del_y$ are $-X_1^{p-1-c}X_2^{p-1-d}\del_y$ and $X_1^{p-1-c}X_2^{p-1-d}\del_x$, so that
$$a_1=C\,\del_y,\ \ \ \ a_2=-C\,\del_x,\ \ \ \ \ C=(x-X_1)^{p-1}(y-X_2)^{p-1}.$$
The $p$-th powers of $x-X_1$ and $y-X_2$ vanish. Any coefficient of the bracket of two differentiated terms in (i) has total degree $4p-5-|\alpha|-|\beta|\ge 2p-1$ in these differences, whereas every nonzero monomial in them has degree at most $2p-2$. This proves (i).

For $F,H\in\O\otimes U$ we have $\D(FH)=(\D F)H+F(\D H)+\sum_i[a_i,F]\del_iH$. If $r<p$ then the coefficients of $\D^r(f)$, $f\in \O$, are sums of products of derivatives $\del^\beta a_i$ with $|\beta|\le r-1$, and these commute with $a_1,a_2$ by part~(i). By induction on $r$ we get $\D^r(fh)=\sum_j\binom{r}{j}\D^j(f)\D^{r-j}(h)$ for $r\le p$. Since the intermediate binomial coefficients vanish for $r=p$, statement~(ii) follows.
\end{proof}
\begin{Lemma}\label{L:3.2}
For every $f\in\O$ we have $\R\big((f\D)^p\big)=f(0)^p\,\R(\D^p)$. Consequently,
\begin{equation}\label{eq:3.2}
	\R\Big(\textstyle{\sum_{j=0}^{p-1}}\,\D^j\varphi\D^{p-j}\Big)\,=\,0\ \ \ \ \ \ (\forall\,\varphi\in \O).
\end{equation}
\end{Lemma}
\begin{proof}
If $\delta$ is a derivation of a commutative ring $R$ of characteristic $p$, then in the Ore extension $R[T;\delta]$, in which $Tr=rT+\delta(r)$ for $r\in R$, one has
$(fT)^p=f^pT^p+(f\delta)^{p-1}(f)T$ for all $f\in R$. To see this, adjoin $t$ with $\delta(t)=1$. The equality ${\rm ad}(fT)^p={\rm ad}\big((fT)^p\big)$ yields $[(fT)^p,t]=(f\delta)^{p-1}(f)$, and comparing with $[T^j,t]=jT^{j-1}$ shows that the Ore normal form of $(fT)^p-f^pT^p$ has no terms of order $2,\ldots,p-1$. Its constant term is zero because both expressions are right divisible by $T$.

Write $(f\D)^p$ and $f^p\D^p+(f\D)^{p-1}(f)\D$ as differential operators with all coefficients placed to the left of $\del_x,\del_y$. In both expressions every coefficient of a term of positive order is a sum of products of derivatives $\del^\beta f$ and $\del^\gamma a_i$ in which the derivatives of $a_1,a_2$ have total order at most $p-1$. By Lemma~\ref{L:3.1}(i) any two such factors commute, so these coefficients are obtained by evaluating the same universal polynomials as in the case where $a_1,a_2$ and their derivatives commute. In that case the identity of the previous paragraph applies to the commutative ring generated by the coefficients, with $\delta=\D$. Hence $(f\D)^p=f^p\D^p+c\D$ with $c=(f\D)^{p-1}(f)$. The derivatives of $a_1,a_2$ occurring in $c$ have total order at most $p-2$, so $c$ commutes with $a_1$ and $a_2$ by Lemma~\ref{L:3.1}(i). It follows that $\R(c\D)={\rm Res}(a_1ca_2-a_2ca_1)=0$. As $f^p=f(0)^p$ in $\O$, the first statement follows. Taking $f=1+\epsilon\varphi$ over $\k[\epsilon]/(\epsilon^2)$ and extracting the coefficient of $\epsilon$ we obtain \eqref{eq:3.2}.
\end{proof}
\begin{Proposition}\label{P:3.3}
The element $u$ is central in $U$ and ${\rm gr}(u)=F_1$.
\end{Proposition}
\begin{proof}
Let $X\in \L$ and $\varphi={\rm div}(X)$. For $\Phi\in\End_\k(\O)\otimes U$ let $V_\Phi=\Phi(x)\otimes\del_x+\Phi(y)\otimes\del_y\in\L\otimes U$, where we use $\O\otimes U\otimes\{\del_x,\del_y\}=\L\otimes U$, and let ${\rm mult}\colon\L\otimes U\to U$ be the multiplication map. Since $\sum_b\omega(g_b,Y)e_b=Y$ for all $Y\in\L$, we have $\R(\Phi)={\rm mult}(V_\Phi)$.

On $\End_\k(\O)\otimes U$ and on $\L\otimes U$ consider the derivations $\Theta_X=[X,\,\cdot\,]\otimes 1+1\otimes \ad X$. Then $[X,{\rm mult}(V)]={\rm mult}(\Theta_X V)$ for $V\in\L\otimes U$. Writing $(x_1,x_2)=(x,y)$, a direct comparison of coefficients gives
$$\Theta_X(V_\Phi)=V_{\Theta_X\Phi}+\textstyle{\sum_i}\Big(\Phi\big(X(x_i)\big)-\sum_j\Phi(x_j)\,\del_j\big(X(x_i)\big)\Big)\otimes\del_i .$$
If $f\mapsto\Phi(f)$ is a derivation $\O\to\O\otimes U$, then $\Phi(X(x_i))=\sum_j\del_j(X(x_i))\Phi(x_j)$ and the last sum vanishes. By Lemma~\ref{L:3.1}(ii) this applies to $\Phi=\D^p$, so that $[X,u]=\R\big(\Theta_X(\D^p)\big)$.

Pairing with $\omega$ in the first factor and using the identity $\omega([X,A],B)+\omega(A,[X,B])=-2\omega(\varphi A,B)$ from the proof of Lemma~\ref{L:2.1}, one finds $\Theta_X(\D)=-2\varphi\D$, where $\varphi$ acts by multiplication. Hence $\Theta_X(\D^p)=-2\sum_{j=0}^{p-1}\D^j\varphi\D^{p-j}$, and $[X,u]=0$ by \eqref{eq:3.2}.

Finally let $\sigma\colon U\to S(\L)$ be the principal symbol and ${\rm ev}_\xi\colon S(\L)\to\k$ evaluation at $\xi\in\L^*$. Applying ${\rm ev}_\xi\circ\sigma$ to the second factor sends $\D$ to $\sum_b\xi(e_b)g_b=J\xi$ and $a_i$ to $(J\xi)(x_i)$. All terms of $\D^p$ have degree $p$ in $U$, so the symbol of $\D^p$ at $\xi$ is the $p$-th power of the derivation $J\xi$, that is, $(J\xi)^{[p]}$. Therefore the symbol of $u$ at $\xi$ is ${\rm Res}\big((J\xi)(x)\,(J\xi)^{[p]}(y)-(J\xi)(y)\,(J\xi)^{[p]}(x)\big)=\omega\big(J\xi,(J\xi)^{[p]}\big)=F_1(\xi)$, and ${\rm gr}(u)=F_1$.
\end{proof}

\subsection{}
\label{ss:3.3}
To lift $F_\psi$ we compare the ordering and symmetrisation maps over $\Z_{(p)}$ up to order $p-2$. Throughout this subsection $R_0=\Z_{(p)}$, and $\l$ denotes a Lie algebra over $R_0$ which is free as an $R_0$-module with a totally ordered basis $\{e_i\,|\,\,i\in I\}$; the index set $I$ may be infinite. We write $[e_i,e_j]=\sum_kc_{ij}^ke_k$.

\smallskip

\noindent
{\it The Rees algebra.} Let $h$ be an indeterminate and let $U_h$ be the quotient of the tensor algebra $T_{R_0}(\l)\otimes_{R_0}R_0[h]$ by the relations $e_je_i-e_ie_j=h[e_j,e_i]$ for all $i,j$. By the PBW theorem, the ordered monomials $e^{\alpha}=\prod_ie_i^{\alpha_i}$ (product taken in the order of $I$) form a basis of $U_h$ over $R_0[h]$, and $U_h/(h-1)U_h\cong U(\l)$. Giving $h$ and all $e_i$ degree one makes $U_h$ a graded algebra; the image of the homogeneous component of degree $\ell$ in $U(\l)$ is the filtration component $U_\ell(\l)$. Let $S=S_{R_0}(\l)=R_0[s_i\,|\,\,i\in I]$, where $s_i$ is the symbol of $e_i$, and let ${\rm ord}_h\colon S[h]\to U_h$ be the $R_0[h]$-linear map sending $s^\alpha$ to $e^\alpha$. It is an isomorphism of graded $R_0[h]$-modules. Let ${\rm ord}\colon S\to U(\l)$ be its specialisation at $h=1$.

\smallskip

\noindent
{\it Comparison over $\Q$.} Put $S_\Q=S\otimes_{R_0}\Q$ and $U_{h,\Q}=U_h\otimes_{R_0}\Q$, and let ${\rm sym}_h\colon S_\Q[h]\to U_{h,\Q}$ be the symmetrisation map, which sends $s_{i_1}\cdots s_{i_\ell}$ to $\frac{1}{\ell!}\sum_{\sigma\in\mathfrak{S}_\ell}e_{i_{\sigma(1)}}\cdots e_{i_{\sigma(\ell)}}$. Let $t=(t_i)_{i\in I}$ be commuting indeterminates and work in the completions with respect to the $t$-adic topology; only finitely many $t_i$ occur in any coefficient we extract. Then ${\rm ord}_h\big(\exp(\sum_it_is_i)\big)=\prod_i\exp(t_ie_i)$ and ${\rm sym}_h\big(\exp(\sum_iy_is_i)\big)=\exp(\sum_iy_ie_i)$ for any family of scalars $y_i$ in the coefficient ring. In $U_{h,\Q}$ the commutator of two elements of $h\l$ equals $h$ times their bracket in $\l$, so the Baker--Campbell--Hausdorff formula gives
\begin{equation}\label{eq:3.bch}
	\log\textstyle{\prod_i}\,\exp(t_ie_i)\,=\,\sum_it_ie_i+\sum_{r\ge 1}h^rC_{r+1}(t),\qquad C_{r+1}(t)=\sum_kC_{r+1}^k(t)\,e_k,
\end{equation}
where $C_{r+1}^k(t)$ is a homogeneous polynomial of degree $r+1$ in the $t_i$, obtained from the degree $r+1$ part of the BCH series by evaluating iterated brackets with the structure constants $c_{ij}^k$. The degree $r+1$ part of the BCH series is a Lie polynomial whose coefficients have denominators dividing $(r+1)!$, so $C^k_{r+1}(t)\in R_0[t]$ for $r+1\le p-1$. Put $b_r(t,s)=\sum_kC^k_{r+1}(t)s_k$; it is linear in $s$. Since the exponent in \eqref{eq:3.bch} is a linear combination of the $e_k$, we obtain
\begin{equation}\label{eq:3.ordsym}
	{\rm ord}_h\big(e^{\sum t_is_i}\big)\,=\,{\rm sym}_h\big(e^{\sum t_is_i}\,E(t,s,h)\big),\qquad E\,=\,\exp\Big(\textstyle{\sum_{r\ge 1}}\,h^rb_r(t,s)\Big).
\end{equation}
Modulo $h^{p-1}$ only the terms $b_1,\ldots,b_{p-2}$ and the powers $(\cdot)^n/n!$ with $n\le p-2$ contribute to $E$, so that
$$E\,\equiv\,\textstyle{\sum_{r=0}^{p-2}}\,h^r\sum_\alpha c_{r,\alpha}(s)\,t^\alpha\pmod{h^{p-1}},\qquad c_{r,\alpha}(s)\in S,$$
where $\alpha$ runs over multi-indices with $|\alpha|\le 2r$ and $c_{r,\alpha}$ is homogeneous of degree $|\alpha|-r$ in $s$. Define differential operators on $S$ by
$$T_r\,=\,\textstyle{\sum_\alpha}\,c_{r,\alpha}(s)\,\del_s^\alpha\ \ (0\le r\le p-2),\qquad Q_0={\rm Id},\ \ Q_r=-\textstyle{\sum_{i=1}^r}T_iQ_{r-i}\ \ (1\le r\le p-2).$$
They are defined over $R_0$, they lower the polynomial degree by $r$, and $T(h)=\sum_rh^rT_r$ and $Q(h)=\sum_rh^rQ_r$ are inverse to each other modulo $h^{p-1}$.
\begin{Lemma}\label{L:3.4a}
For every $F\in S_\Q$ one has ${\rm ord}_h(F)\equiv{\rm sym}_h\big(T(h)F\big)$ and ${\rm sym}_h(F)\equiv{\rm ord}_h\big(Q(h)F\big)$ modulo $h^{p-1}U_{h,\Q}$.
\end{Lemma}
\begin{proof}
Since $t^\alpha e^{\sum t_is_i}=\del_s^\alpha e^{\sum t_is_i}$, the right hand side of \eqref{eq:3.ordsym} equals ${\rm sym}_h\big(T(h)e^{\sum t_is_i}\big)$ modulo $h^{p-1}$. Both sides of the resulting identity are linear in $e^{\sum t_is_i}=\sum_\beta t^\beta s^\beta/\beta!$, and $T(h)$ does not involve $t$. Comparing the coefficients of $t^\beta/\beta!$ we get the first congruence for $F=s^\beta$, hence for all $F$. Applying it to $Q(h)F$ gives the second.
\end{proof}
\begin{Lemma}\label{L:3.4}
Let $F\in S$ be homogeneous of degree $\ell$ and $Y\in\l$, and put $z(F)={\rm ord}\big(\sum_{r=0}^{p-2}Q_rF\big)\in U(\l)$. Then
$$[Y,z(F)]-z(Y\cdot F)\in U_{\ell-p+1}(\l),$$
where $Y\cdot F$ denotes the adjoint action of $Y$ on $S$.
\end{Lemma}
\begin{proof}
Let $\delta_Y$ be the $R_0[h]$-linear derivation of $U_h$ with $\delta_Y(e_i)=[Y,e_i]$; it is well defined because it preserves the defining relations. Then $[Y,u]=h\,\delta_Y(u)$ for all $u\in U_h$, and $\delta_Y$ specialises to $\ad Y$ at $h=1$. Since ${\rm sym}_h$ is the average over all orderings and $\delta_Y$ is a derivation, $\delta_Y\big({\rm sym}_h(F)\big)={\rm sym}_h(Y\cdot F)$ for all $F\in S_\Q$.

Put $z_h(F)={\rm ord}_h\big(Q(h)F\big)\in U_h$, a homogeneous element of degree $\ell$. By Lemma~\ref{L:3.4a}, $z_h(F)={\rm sym}_h(F)+h^{p-1}x$ and $z_h(Y\cdot F)={\rm sym}_h(Y\cdot F)+h^{p-1}x'$ with $x,x'\in U_{h,\Q}$. Hence
$$\delta_Y\big(z_h(F)\big)-z_h(Y\cdot F)\,=\,h^{p-1}\big(\delta_Y(x)-x'\big).$$
The left hand side lies in $U_h$. Its coordinates in the ordered PBW basis are polynomials in $h$ with coefficients in $R_0$ that are divisible by $h^{p-1}$ in $\Q[h]$, hence in $R_0[h]$. Therefore $\delta_Y(z_h(F))-z_h(Y\cdot F)=h^{p-1}a$ with $a\in U_h$ homogeneous of degree $\ell-(p-1)$. Specialising at $h=1$ gives the lemma.
\end{proof}

\smallskip

\noindent
{\it Descent to $\L$.} Let $\l=R_0[x,y]\del_x\oplus R_0[x,y]\del_y$ with the monomial basis $x^iy^j\del_x$, $x^iy^j\del_y$, in any fixed order. The $R_0$-submodule $\mathfrak{i}=p\l+(x^p,y^p)\l$ is an ideal of $\l$. Indeed, $[x^pf\del_a,g\del_b]=x^pf\,\del_a(g)\del_b-g\,\del_b(x^pf)\del_a$, and $\del_b(x^pf)\in x^pR_0[x,y]+pR_0[x,y]$; the same holds with $y^p$ in place of $x^p$. Moreover $\l/\mathfrak{i}\cong\L$. The structure constants of $\l$ and those of $\L$ agree modulo $\mathfrak{i}$, so the operators $T_r,Q_r$ preserve the ideal of $S$ generated by $\mathfrak{i}$, and ${\rm ord}$ maps it into the two-sided ideal of $U(\l)$ generated by $\mathfrak{i}$. Hence $T_r$, $Q_r$, $z$ and Lemma~\ref{L:3.4} descend to $S(\L)$ and $U(\L)$. All of them preserve the integer grading of \ref{ss:3.1}, since the structure constants do.

\smallskip

\noindent
{\it The lift of $F_\psi$.} Put
$$v'\,:=\,z(F_\psi)\,=\,{\rm ord}\big(\textstyle{\sum_{r=0}^{p-2}}\,Q_rF_\psi\big)\in U.$$
Since $Q_r$ lowers the degree by $r$, ${\rm gr}(v')=F_\psi$, and $v'$ is homogeneous of integer degree $\theta_2$. Let $Y\in\L$ be homogeneous of integer degree $d\ge -1$. As $F_\psi$ is invariant, Lemma~\ref{L:3.4} gives $[Y,v']\in U_{(p-1)^2}$, because $\ell-p+1=(p-1)^2$ for $\ell=p^2-p$. The element $[Y,v']$ is homogeneous of integer degree $\theta_2+d\ge \theta_2-1$. On the other hand, every element of $\L$ has integer degree at most $2p-3$, so every element of $U_{(p-1)^2}$ has integer degree at most $(p-1)^2(2p-3)$. Since
$$\theta_2-1-(p-1)^2(2p-3)\,=\,2p(p-1)(p-2)-1-(p-1)^2(2p-3)\,=\,p^2-4p+2\,>\,0,$$
we conclude that $[Y,v']=0$. Thus $v'$ is central. Finally, we set
$$v\,:=\,v'+h_0(\underline{z}).$$
This is a central element of $U$ with ${\rm gr}(v)=F_\psi+h_0^p=F_2$.

\subsection{}
\label{ss:3.4}
Since principal symbols of central elements are $\L$-invariant, ${\rm gr}\,Z\subseteq S(\L)^\L$. On the other hand, ${\rm gr}\,Z$ contains $S(\L)^p={\rm gr}\,Z_p(\L)$ as well as $F_1={\rm gr}(u)$ and $F_2={\rm gr}(v)$. As $S(\L)^\L=S(\L)^p[F_1,F_2]$, it follows that ${\rm gr}\,Z=S(\L)^\L$. If $c\in Z$ has filtration degree $\ell$, write ${\rm gr}(c)=\sum_{i,j}a_{ij}F_1^iF_2^j$ with $a_{ij}\in S(\L)^p$ and $0\le i,j\le p-1$, lift each $a_{ij}$ to $Z_p(\L)$, and subtract the corresponding $Z_p(\L)$-linear combination of the $u^iv^j$. Induction on $\ell$ shows that
$$Z\,=\,\textstyle{\bigoplus_{0\le i,j\le p-1}}\,Z_p(\L)\,u^iv^j.$$
The sum is direct since a nontrivial relation would contradict, in its highest PBW degree, the $S(\L)^p$-linear independence of the $F_1^iF_2^j$. The same argument shows that the basis is strict for the PBW filtration:
\begin{equation}\label{eq:3.3}
	\deg\textstyle{\sum_{i,j}}\,c_{ij}u^iv^j\,=\,\max_{c_{ij}\ne 0}\big(p\deg_z(c_{ij})+i(p+1)+j(p^2-p)\big),
\end{equation}
where $\deg_z$ denotes the degree of an element of $Z_p(\L)$ as a polynomial in the $z_b$. This proves the first two assertions of Theorem~\ref{main1}; the semi-invariance of $u$ and $v$ is established in \ref{ss:3.7}.

\subsection{}
\label{ss:3.5}
We need a description of the leading term of $a^p$ for $a\in U$. We use the restricted Poisson structure on $S(\L)$ attached to the PBW filtration; see \cite[Theorems~6.1, 8.2 and Lemma~7.2]{BYZ}.
\begin{Lemma}\label{L:3.5}
Let $f\in S^\ell(\L)$ have coefficients in $\mathbb{F}_p$ and let $a\in U_\ell$ have principal symbol $f$. Then $a^p-f(\underline{z})\in U_{p\ell-p+1}$. If $f\in S(\L)^\L$ then the image $f^{\{p\}}$ of $a^p-f(\underline{z})$ in $S^{p\ell-p+1}(\L)$ satisfies
\begin{equation}\label{eq:3.4}
	\big(df^{\{p\}}\big)_\xi\,=\,\big((df)_\xi\big)^{[p]}\ \ \ \ \ \ (\forall\,\xi\in\L^*).
\end{equation}
\end{Lemma}
\begin{proof}
Let $\mathcal{R}_h=\bigoplus_dh^dU_d$ be the Rees algebra of the PBW filtration, write $\tilde x=hx$ for $x\in\L$, and let $a_h=h^\ell a\in\mathcal{R}_h$. The elements $\tilde x^p-h^{p-1}\widetilde{x^{[p]}}=h^p(x^p-x^{[p]})$ are central with symbol $x^p$. In $\mathcal{R}_h$ every commutator carries a factor $h$; modulo $h^{p-1}$ all $p$-fold commutators and all $p$-th powers of commutators vanish, so that the $p$-th power map is additive and multiplicative in this quotient by \cite[Lemma~7.2]{BYZ}. Hence $a_h^p-h^{p\ell}f(\underline{z})$ is divisible by $h^{p-1}$. Its coefficient of $h^{p-1}$ does not change if we replace $a_h$ by $a_h+hb$: by Jacobson's formula \cite{Jac} the difference $(a_h+hb)^p-a_h^p-(hb)^p$ is a sum of commutators of length $p$, each containing at least one factor $hb$ and $p-1$ brackets, hence divisible by $h^p$. This gives the first statement. The resulting operation $f\mapsto f^{\{p\}}$ extends the $p$-mapping of $\L$ \cite[Theorem~6.1]{BYZ}, and \cite[Theorem~8.2]{BYZ} yields $(df)^{[p]}=d(f^{\{p\}})$ on $\Omega_{S(\L)/\k}=S(\L)\otimes\L$. At $\xi$ the isotropy quotient of this restricted Lie--Rinehart algebra\footnote{Here $\Omega_{S(\L)/\k}=S(\L)\otimes\L$ carries the bracket extending that of $\L$ and the anchor $1\otimes x\mapsto\{x,\cdot\,\}$; its isotropy at $\xi$ is the set of elements whose anchor vanishes at $\xi$, modulo $\m_\xi\Omega_{S(\L)/\k}$, where $\m_\xi\subset S(\L)$ is the maximal ideal of $\xi$. See \cite[Sections~6--8]{BYZ}.} is $\L_\xi$ with its $p$-mapping; indeed, the formula $(aD)^{[p]}=a^pD^{[p]}+(aD)^{p-1}(a)D$ shows that the submodule vanishing at $\xi$ is preserved whenever the anchor vanishes at $\xi$. An invariant differential lies in the isotropy. This implies \eqref{eq:3.4}.
\end{proof}
Let $\iota$ be the principal anti-automorphism of $U$, so that $\iota(x)=-x$ for $x\in\L$. Then $\iota(z_b)=-z_b$ for all $b$. Since $\iota$ is an anti-automorphism, it preserves $Z$, and it acts on principal symbols of degree $d$ by $(-1)^d$. Replacing $u$ by $\frac12\big(u+\iota(u)\big)$ and $v'$ by $\frac12\big(v'+\iota(v')\big)$, which does not change their principal symbols $F_1$ and $F_\psi$ (as $p+1$ and $p^2-p$ are even) nor their centrality and homogeneity, we may and shall assume that $\iota(u)=u$ and $\iota(v')=v'$. Since $h_0$ has even degree $p-1$, also $\iota(h_0(\underline{z}))=h_0(\underline{z})$, and hence $\iota(v)=v$. In particular, all coefficients in the expansions of $u^p$ and $v^p$ with respect to the basis \eqref{eq:3.3} are even polynomials in the $z_b$.

\subsection{}
\label{ss:3.6}
We now prove the first relation in \eqref{eq:1}.
\begin{Proposition}\label{P:3.6}
One has $u^p+uv=F_1(\underline{z})$.
\end{Proposition}
\begin{proof}
By Lemma~\ref{L:3.5}, $u^p-F_1(\underline{z})$ has PBW degree at most $p^2+1$, and it is homogeneous of integer degree $p\theta_1=2p^2(p-2)$. Expand it as $\sum_{i,j}c_{ij}u^iv^j$ with respect to \eqref{eq:3.3}, and let $k_{ij}=\deg_z(c_{ij})$. Then $c_{ij}$ has integer degree $(p-i)\theta_1-j\theta_2$, while each $z_b$ has integer degree at most $p(2p-3)$; moreover $pk_{ij}+i(p+1)+j(p^2-p)\le p^2+1$. By \ref{ss:3.5}, $u^p-F_1(\underline{z})$ is fixed by $\iota$, so $k_{ij}$ is even whenever $c_{ij}\ne0$.
\begin{itemize}
\item[(a)] If $j=1$, the PBW bound gives $(k_{ij},i)\in\{(0,0),(1,0),(0,1)\}$. The case $(0,0)$ has the wrong integer degree, and $(1,0)$ has odd $z$-degree. Only $c_{11}=:\kappa\in\k$ remains.
\item[(b)] If $j=0$ and $i=0$, then $k_{00}\le p$ and $k_{00}\,p(2p-3)\ge 2p^2(p-2)$ force $k_{00}=p$, which is odd; so $c_{00}=0$.
\item[(c)] If $j=0$ and $i=1$, the same two inequalities give $k_{10}=p-1$; put $B=c_{10}$.
\item[(d)] If $j=0$ and $i\ge2$, the PBW bound gives $k_{i0}\le p-1-i$, whereas the integer degree requires $k_{i0}\ge2(p-i)(p-2)/(2p-3)>p-1-i$; so $c_{i0}=0$.
\end{itemize}
Hence
$$u^p\,=\,F_1(\underline{z})+Bu+\kappa uv,\ \ \ \ \ \ B\in Z_p(\L),\ \ \deg_z(B)= p-1,\ \ \kappa\in\k.$$

Let $\xi\in\L^*$, $X=J\xi$, and put $h_1=(dF_1)_\xi$ and $h_2=(dF_\psi)_\xi$, elements of $\L$. Let $\zeta\in\L^*$. As in the proof of Proposition~\ref{P:2.2}, $\zeta(h_1)=\omega(J\zeta,X^{[p]})+\omega\big(X,(\ad X)^{p-1}J\zeta\big)$. The identity in the proof of Lemma~\ref{L:2.1} gives $\omega\big(Y_1,[X,Y_2]\big)=-\omega\big((\ad X+2{\rm div}(X))Y_1,Y_2\big)$ for $Y_1,Y_2\in\L$, and $(\ad X+2{\rm div}(X))(aX)=\big(X(a)+2{\rm div}(X)a\big)X$ for $a\in\O$. Hence $\omega\big(X,(\ad X)^{p-1}J\zeta\big)=(-1)^{p-1}\omega\big((\ad X+2{\rm div}(X))^{p-1}X,J\zeta\big)=\omega(q_XX,J\zeta)$, and therefore $h_1=X^{[p]}-q_XX$. By \eqref{eq:2.4}, ${\rm Res}\big(a\,\nu(X,h_2)\big)=0$ for all $a\in\O$, whence $\nu(X,h_2)=0$. On the nonempty open set where $w\ne 0$ the pair $X,X^{[p]}$ is an $\O$-basis of $\L$, so $h_2\in\O X$. The characteristic polynomial of $X$ is an additive polynomial. Indeed, if $X$ has a simple spectrum consisting of units, its eigenvalues form a two-dimensional $\mathbb{F}_p$-subspace of $\k$, whose root polynomial is additive; such derivations form a nonempty open subset of $\L$ containing the elements used in the proof of \eqref{eq:2.4}, and the coefficients of $\det(T\,{\rm Id}-X)$ are polynomial in $X$. By the Cayley--Hamilton theorem,
$X^{[p^2]}+F_\psi(\xi)X^{[p]}+F_0(\xi)X=0$, where $F_0(\xi)$ is the coefficient of $T$ in $\det(T\,{\rm Id}-X)$. Since $\O X$ is an ideal of the Lie algebra $\k X^{[p]}+\O X$, Jacobson's formula \cite{Jac} gives $(X^{[p]}-q_XX)^{[p]}\equiv X^{[p^2]}-(q_XX)^{[p]}\pmod{\O X}$, and $(q_XX)^{[p]}\equiv q_X(0)^pX^{[p]}=h_0(\xi)^pX^{[p]}\pmod{\O X}$ by Hochschild's formula. Therefore
$$h_1^{[p]}\,\equiv\,-\big(F_\psi(\xi)+h_0(\xi)^p\big)h_1\,=\,-F_2(\xi)h_1\pmod{\O X}.$$
On the other hand, applying \eqref{eq:3.4} to the correction term $B_{\rm sym}F_1+\kappa F_1F_2$ we get $h_1^{[p]}=(B_{\rm sym}+\kappa F_2)h_1+\kappa F_1h_2$ with $B_{\rm sym}\in S(\L)^p$. Comparing the coefficients of $X^{[p]}$ yields $B_{\rm sym}+(\kappa+1)F_2=0$. Since $(dF_2)_\eta\ne 0$, we have $F_2\notin S(\L)^p$, forcing $\kappa=-1$ and $B_{\rm sym}=0$. The integer-degree bound excludes lower layers of $B$, so $B=0$.
\end{proof}

\subsection{}
\label{ss:3.7}
We now turn to the action of $G$. For $g\in G$ write $\delta(g)=\det(dg_0)$.
\begin{Lemma}\label{L:3.7}
For all $g\in G$ we have $g(u)=\delta(g)^{p(p-2)}u$ and $g(v)=\delta(g)^{p(p-1)(p-2)}v$. Moreover, $w$ is a $G$-semi-invariant of weight $p^2-p-1$.
\end{Lemma}
\begin{proof}
The degree $p+1$ component of $S(\L)^\L$ is $\k F_1$ by Theorem~\ref{T:2.6}, and $g(F_1)=\delta(g)^{p(p-2)}F_1$ by Lemma~\ref{L:2.4}. Hence $g(u)-\delta(g)^{p(p-2)}u$ is a central element of PBW degree at most $p$ fixed by $\iota$. Such an element is a scalar plus a linear combination of the $z_b$. Since $\iota$ negates the $z_b$, it is a scalar, and since $u$ lies in the augmentation ideal of $U$, this scalar is zero. As $G$ acts on $Z_p(\L)$ compatibly with its action on $S(\L)$, and $U$ is a domain, cancelling $u$ in $g(u)^p+g(u)g(v)=g\big(F_1(\underline{z})\big)$ gives the formula for $g(v)$. The last assertion is Lemma~\ref{L:2.4}.
\end{proof}
This completes the proof of Theorem~\ref{main1}. Since every $G$-semi-invariant central element is a polynomial in $u$ and $v$ by induction on the filtration degree, $Z_{\rm HC}(U)=\k[u,v]$.

Next we consider the affine plane $\{\xi_{a,b}\,|\,\,a,b\in\k\}\subset\L^*$, where $\xi_{a,b}=\eta+a\,\del_y^*+b\,\del_x^*$; we call it the {\it slice}. Its residue dual and characteristic polynomial are
\begin{equation}\label{eq:3.5}
	J\xi_{a,b}=\Delta+ax^{p-1}y^{p-1}\del_x-bx^{p-1}y^{p-1}\del_y,\ \ \ \ \ \det(T\,{\rm Id}-J\xi_{a,b})=T^{p^2}+bT^p-aT.
\end{equation}
To see this, put $\tau=x^{p-1}y^{p-1}$ and $c(h)=\big(a\del_x(h)-b\del_y(h)\big)(0)$ for $h\in\O$. As $\tau h=h(0)\tau$, we have $J\xi_{a,b}=\Delta+\rho$, where $\rho(h)=c(h)\,\tau$ is an operator of rank one. Hence, for every $n\ge1$,
$$(J\xi_{a,b})^n(h)=\Delta^n(h)+\textstyle{\sum_{k=0}^{n-1}}\,c\big((J\xi_{a,b})^k(h)\big)\,\Delta^{n-1-k}(\tau).$$
By \eqref{eq:2.5}, $\Delta$ is nilpotent with $\Delta^{p^2-2}(\tau)=x$ and $\Delta^{p^2-p-1}(\tau)=y$, while $\Delta^k(\tau)$ has no linear terms for all other $k\le p^2-1$; thus $c(\Delta^k(\tau))$ equals $a$ for $k=p^2-2$, $-b$ for $k=p^2-p-1$, and $0$ otherwise. The rank-one determinant identity gives
$$\det(T\,{\rm Id}-J\xi_{a,b})=T^{p^2}\Big(1-\textstyle{\sum_{k\ge0}}\,T^{-k-1}c\big(\Delta^k(\tau)\big)\Big)=T^{p^2}+bT^p-aT.$$
For the $p$-th power take $n=p$ and $h\in\{x,y\}$. For $k\le p-1$ the element $\Delta^k(\tau)$ has degree at least $p-1\ge2$, so $c\big((J\xi_{a,b})^k(h)\big)=c\big(\Delta^k(h)\big)$. Now $c(x)=a$, $\Delta(x)=1$; and $c(y)=-b$, $\Delta^k(y)$ has degree at least two for $1\le k\le p-2$, and $\Delta^{p-1}(y)=-\Delta^{p-2}(x^{p-1})=-(p-1)!\,x=x$ by Wilson's theorem. As $\Delta^{p-1}(\tau)=(p-1)!\,y^{p-1}=-y^{p-1}$ and $\Delta^p=\del_y$, we obtain $(J\xi_{a,b})^p(x)=a\Delta^{p-1}(\tau)=-ay^{p-1}$ and $(J\xi_{a,b})^p(y)=1-b\Delta^{p-1}(\tau)+a\tau$. Since $(J\xi_{a,b})^{[p]}$ is the derivation $(J\xi_{a,b})^p$,
$$(J\xi_{a,b})^{[p]}\,=\,-ay^{p-1}\del_x+\big(1+ax^{p-1}y^{p-1}+by^{p-1}\big)\del_y.$$ Thus $F_1=a$, $F_\psi=b$, $w=1$ and $h_0=0$ on the slice; in particular, $F_2=b$ and $w(\eta)=1$. The map $G\times\A^2\to\L^*$, $(g,a,b)\mapsto g\cdot\xi_{a,b}$, is dominant: its differential at $(1,0,0)$ has image $\L_{(0)}\cdot\eta+\k\del_x^*+\k\del_y^*$, which is all of $\L^*$ because $\L_{(0)}\cdot\eta$ is the annihilator of $\L_\eta=\k\Delta\oplus\k\del_y$ in $\L^*$ (as $\L_\eta\cap\L_{(0)}=0$) and $\del_x^*,\del_y^*$ pair nondegenerately with $\Delta,\del_y$.

Let $f\in S^{q-s}(\L)$, where $s\ge 1$ and $0\le q\le p^2-p$, be a $G$-semi-invariant of weight $q(p-2)$. Under the cocharacter $\mu$ of \ref{ss:1.4} the two nonzero coordinates of $\eta$ have $\mu$-weight $p^2-p-1$, while the two slice coordinates have weights $-p$ and $-1$. For a monomial surviving on the slice let $A,B$ be its exponents in the slice coordinates. Comparing total and $\mu$-weights gives
\begin{equation}\label{eq:3.6}
	s(p^2-p-1)+(p^2-1)A+(p^2-p)B\,=\,q.
\end{equation}
The only solution is $s=1$, $q=p^2-p-1$, $A=B=0$. Consequently, $f=0$ unless $\deg f=p^2-p-2$ and $q=p^2-p-1$, in which case $f$ is constant on the slice. By dominance, two such semi-invariants agreeing on the slice coincide.

\subsection{}
\label{ss:3.8}
We are finally ready to prove our second main result.
\begin{Proposition}\label{P:3.8}
One has $v^p-w(\underline{z})^{p-2}u=F_2(\underline{z})$.
\end{Proposition}
\begin{proof}
Expand $v^p-F_2(\underline{z})=\sum K_{ij}u^iv^j$ with respect to \eqref{eq:3.3}. By Lemma~\ref{L:3.7}, each $K_{ij}$ is a semi-invariant of weight $q(p-2)p$ with $q=p(p-1)-i-(p-1)j$, and $\deg_z(K_{ij})\le q-1$ by Lemma~\ref{L:3.5}. Write $K_{ij}=f_{ij}(\underline{z})$ with $f_{ij}\in S(\L)$; since $g(f(\underline{z}))=(g f)(\underline{z})$ and the weight of $f(\underline{z})$ is $p$ times that of $f$, each homogeneous component of $f_{ij}$ is a $G$-semi-invariant of weight $q(p-2)$ and degree at most $q-1$. Applying \eqref{eq:3.6} we find that only $(i,j)=(1,0)$ survives, so that $v^p=F_2(\underline{z})+K(\underline{z})u$ with $\deg K=p^2-p-2$. At $\eta$ the gradient of $F_2$ is $\Delta$, and \eqref{eq:3.4} applied to the correction term gives $K(\eta)\,\del_y=\Delta^{[p]}=\del_y$, so that $K(\eta)=1$. The polynomial $K$ has the same degree and $G$-weight as $w^{p-2}$, and both restrict to $1$ on the slice. Hence $K=w^{p-2}$.
\end{proof}
\begin{proof}[Proof of Theorem~\ref{main2}]
Propositions~\ref{P:3.6} and \ref{P:3.8} give the relations \eqref{eq:1}. Give the formal variables $U_1,U_2$ degrees $p+1$ and $p^2-p$ and take $Z_p(\L)$ as the coefficient ring. The leading monomials of the two relations are $U_1^p$ and $U_2^p$; they are monic and coprime, so the quotient of $Z_p(\L)[U_1,U_2]$ by the two relations is free over $Z_p(\L)$ with basis $U_1^iU_2^j$, $0\le i,j\le p-1$. The natural map from this quotient to $Z$ sends this basis to the basis of \ref{ss:3.4} and hence is an isomorphism. The relations form a regular sequence, so ${\rm Specm}\,Z$ is a complete intersection in $\A^{n+2}(\k)$.
\end{proof}

\subsection{}
\label{ss:3.9}
At $\eta$ we have $F_1(\eta)=F_2(\eta)=0$ and $w(\eta)=1$. Hence the fibre of $Z$ over $\eta$ is
$$Z_\eta\,=\,\k[u,v]/(u^p+uv,\,v^p-u)\,\cong\,\k[X]/\big(X^{p^2}+X^{p+1}\big).$$
\begin{Proposition}\label{P:3.9} Let $N=p^{p^2-1}$. Then
$$U_\eta(\L)\,\cong\,{\rm Mat}_N(Z_\eta)\,\cong\,{\rm Mat}_N\big(\k[X]/(X^{p+1})\big)\oplus\,{\rm Mat}_N(\k)^{\oplus(p^2-p-1)}.$$
\end{Proposition}
\begin{proof}
We follow the strategy of \cite[2.1]{PS99}. Let $t$ be an indeterminate, let $\mathcal{R}=\bigoplus_{d\ge0}t^dU_d\subset U[t]$ be the Rees algebra of the PBW filtration, and write $\tilde x=tx$ for $x\in\L$, so that $\tilde x\tilde y-\tilde y\tilde x=t\,\widetilde{[x,y]}$ in $\mathcal{R}$. Put
$$\mathcal{A}\,=\,\mathcal{R}\big/\big(\tilde x^p-t^{p-1}\widetilde{x^{[p]}}-\eta(x)^p\,\big|\,\,x\in\L\big),\qquad \mathcal{A}_c=\mathcal{A}/(t-c)\mathcal{A}\ \ (c\in\k).$$
By the PBW theorem, $\mathcal{A}$ is a free $\k[t]$-module with basis the ordered monomials $\tilde e^\alpha$, $0\le\alpha_i\le p-1$, and the $\k[t]$-subalgebra $\mathcal{A}_{(0)}$ generated by the $\tilde x$ with $x\in\L_{(0)}$ is free with basis the ordered monomials in a basis of $\L_{(0)}$. We have $\mathcal{A}_1=U_\eta(\L)$ and $\mathcal{A}_0=S_\eta(\L)=U_{\eta,0}(\L)$. For $c\ne0$ the assignment $x\mapsto c^{-1}\tilde x$ induces an isomorphism $U_{c^{-1}\eta}(\L)\cong\mathcal{A}_c$ mapping $U_{c^{-1}\eta}(\L_{(0)})$ onto the image of $\mathcal{A}_{(0)}$. Choosing $s\in\k^\times$ with $s^{-(p^2-p-1)}=c^{-1}$ and applying $\mu(s)$, which preserves $\L_{(0)}$ and maps $\eta$ to $c^{-1}\eta$, we obtain an isomorphism $U_\eta(\L)\cong\mathcal{A}_c$ preserving the subalgebras generated by $\L_{(0)}$.

The elements $\tilde u=t^{p+1}u$ and $\tilde v=t^{p^2-p}v$ are central in $\mathcal{R}$, hence in $\mathcal{A}$. Their images in $\mathcal{A}_1$ are $u$ and $v$, their images in $\mathcal{A}_0$ are ${}^\eta F_1$ and ${}^\eta F_2$, and by Lemma~\ref{L:3.7} the isomorphisms of the previous paragraph map them to nonzero multiples of $u$ and $v$. Consider the $\k[t]$-linear map
$$m\colon\ \textstyle{\bigoplus_{0\le i,j\le p-1}}\,\tilde u^i\tilde v^j\,\mathcal{A}_{(0)}\longrightarrow\mathcal{A},$$
where the source is the free $\k[t]$-module $\k[t]^{p^2}\otimes_{\k[t]}\mathcal{A}_{(0)}$ and $m$ is multiplication. Source and target are free of rank $p^{2p^2}$. At $t=0$, $m$ becomes the multiplication map $\k[{}^\eta F_1,{}^\eta F_2]\otimes S_\eta(\L_{(0)})\to S_\eta(\L)$. This map is an isomorphism: both sides have dimension $p^{2p^2}$, and it is surjective because $S_\eta(\L)$ is a local algebra whose cotangent space is spanned by the images of $\L_{(0)}$ together with $(dF_1)_\eta=\del_y$ and $(dF_2)_\eta=\Delta$, by Lemma~\ref{L:2.3} and $\L=\L_\eta\oplus\L_{(0)}$. Hence $\det m\in\k[t]$ does not vanish at $t=0$, so $m$ is an isomorphism at some $c\ne0$, and by the isomorphisms above at $c=1$. Therefore the images of $u^iv^j$, $0\le i,j\le p-1$, are linearly independent in $U_\eta(\L)$, they span $Z_\eta$, and
$$U_\eta(\L)\,\cong\,Z_\eta\otimes U_\eta(\L_{(0)}).$$
Since $\L=\L_\eta\oplus\L_{(0)}$, the stabiliser of $\eta_{\vert\L_{(0)}}$ in $\L_{(0)}$ is trivial, and \cite[Theorem~4.4(3)]{PS99} gives $U_\eta(\L_{(0)})\cong{\rm Mat}_N(\k)$. Thus $U_\eta(\L)\cong{\rm Mat}_N(Z_\eta)$. Since $X^{p^2}+X^{p+1}=X^{p+1}(X^{p^2-p-1}+1)$ and $p\nmid p^2-p-1$, the Chinese Remainder Theorem completes the proof.
\end{proof}
\noindent
Alternatively, one can deduce that $U_\eta(\L)\cong{\rm Mat}_N(Z_\eta)$ by using Theorem~\ref{main3}. Indeed, since $\eta$ is regular, $U_\eta(\L)$ is an Azumaya algebra of degree $N$ over the Artinian algebra $Z_\eta$, and an Azumaya algebra over an Artinian local ring is a matrix algebra.
If $\zeta$ is a primitive $(p^2-p-1)$-st root of unity then $\mu(\zeta)$ fixes $\eta$, and Lemma~\ref{L:3.7} shows that $\mu(\zeta)$ acts on $Z_\eta$ by $v\mapsto\zeta^{\pm 1} v$. It follows that $\mu(\zeta)$ fixes the local block and permutes the $p^2-p-1$ simple blocks of $U_\eta(\L)$ transitively.

\section{The Azumaya locus of $U(\L)$}
\label{s:4}
\subsection{}
\label{ss:4.1}
Put $d_0=p^2-1$, so that $N=p^{d_0}$. Following \cite[5.2]{PS99} we let $\mathcal{X}_{d_0}$ be the set of all $\chi\in\L^*$ for which $U_\chi(\L)$ has a finite dimensional module of dimension not divisible by $p^{d_0}$. Since all simple $U(\L)$-modules have dimension at most $N$, the set $\mathcal{X}_{d_0}$ coincides with the set $\mathcal{X}_{\rm sing}$ of \ref{ss:1.7}. By \cite[Proposition~5.2]{PS99}, $\mathcal{X}_{d_0}$ is Zariski closed and its associated cone is contained in $\L^*_{\rm sing}$. If $\chi\notin\mathcal{X}_{d_0}$ then every point of $\Psi^{-1}(\chi)$ lies in $\Omega$; see \cite[Theorems~5 and 6]{Zass} and \cite[5.2]{PS99}. We are going to show that $\mathcal{X}_{d_0}\subseteq\L^*_{\rm sing}$.
\begin{Lemma}\label{L:4.1}
Let $t\mapsto g(t)$ be a morphism $\k^\times\to G$, let $r>0$, and suppose that the morphism $\k^\times\to\L^*$, $t\mapsto t^rg(t)\cdot\chi$, extends to a morphism $\A^1\to\L^*$ whose value $\chi_\infty$ at $t=0$ is a nonzero element of $\L^*_{\rm reg}$. Then $\chi\notin\mathcal{X}_{d_0}$.
\end{Lemma}
\begin{proof}
The set $\mathcal{X}_{d_0}$ is $G$-stable. If $\chi\in\mathcal{X}_{d_0}$ then the curve $(g(t)\cdot\chi:1)=(t^rg(t)\cdot\chi:t^r)$ lies in the projective closure of $\mathcal{X}_{d_0}$ in $\mathbb{P}(\L^*\oplus\k)$. Its value $(\chi_\infty:0)$ at $t=0$ then lies in the projectivised associated cone of $\mathcal{X}_{d_0}$ \cite[5.1]{PS99}, which is contained in $\mathbb{P}(\L^*_{\rm sing})$. This is a contradiction.
\end{proof}
For $\alpha_t\in G$ with $\alpha_t(x)=t^ax$ and $\alpha_t(y)=t^by$ we have $j_{\alpha_t}=t^{a+b}=\det(d(\alpha_t)_0)$, so \eqref{eq:2.6} gives
\begin{equation}\label{eq:4.1}
	J(\alpha_t\cdot\chi)\,=\,t^{-(p-2)(a+b)}(\alpha_t)_*(J\chi).
\end{equation}

\subsection{}
\label{ss:4.2}
We first consider the linear functions $\chi$ with $J\chi\notin\L_{(0)}$, using the normal form of Lemma~\ref{L:4.2}.

Let $\chi\in\L^*_{\rm reg}$ with $J\chi\notin\L_{(0)}$, and write $J\chi$ in the form \eqref{eq:4.2}. If $c(0)\ne 0$, take $(a,b)=(1,p)$ in \eqref{eq:4.1}; then
$$t^{p^2-p-1}J(\alpha_t\cdot\chi)\,=\,f(tx,t^py)\big(\del_x+x^{p-1}c(t^py)\del_y\big)\,\longrightarrow\, f(0)\big(\del_x+c(0)x^{p-1}\del_y\big)$$
as $t\to 0$, and the limit is regular by Lemma~\ref{L:4.2}. If $c(0)=0\ne c'(0)$, a change of the coordinate $y$ turns $c(y)\del_y$ into $c'(0)y\del_y$ (the coefficients are determined successively, the relevant diagonal entries $c'(0)(j-1)$, $2\le j\le p-1$, being nonzero), and $(a,b)=(0,1)$ gives the regular limit $f(x,0)\big(\del_x+c'(0)x^{p-1}y\del_y\big)$ after multiplication by $t^{p-2}$. In both cases Lemma~\ref{L:4.1} shows that $\chi\notin\mathcal{X}_{d_0}$.

\subsection{}
\label{ss:4.3}
Next suppose $X=J\chi\in\L_{(0)}$ and let $\lambda_1,\lambda_2$ be the eigenvalues of its linear part $X_0\in\L_0\cong\gl(2,\k)$.
\begin{Lemma}\label{L:4.3}
The linear function $\chi$ is regular if and only if $\lambda_1,\lambda_2$ are linearly independent over $\mathbb{F}_p$. If they are dependent then $(dF_\psi)_\chi=0$.
\end{Lemma}
\begin{proof}
Suppose first that $\lambda_1,\lambda_2$ are independent, so that $X_0$ is diagonalisable; after conjugation by $G_0$ we may take $X_0=\lambda_1x\del_x+\lambda_2y\del_y$. By Lemma~\ref{L:2.1}, $Y\in\L_{J^{-1}X_0}$ if and only if $[Y,X_0]+2{\rm div}(Y)X_0=0$. This equation respects the decomposition of $\L$ into $\t_0$-weight spaces. For $\beta=(\beta_1,\beta_2)$ the weight space $\L^\beta$ is spanned by those of $x^{\beta_1+1}y^{\beta_2}\del_x$ and $x^{\beta_1}y^{\beta_2+1}\del_y$ whose exponents lie in $\{0,\ldots,p-1\}$. Put $\lambda=\lambda_1\beta_1+\lambda_2\beta_2$. For $Y=A\,x^{\beta_1+1}y^{\beta_2}\del_x+B\,x^{\beta_1}y^{\beta_2+1}\del_y$ one has $[Y,X_0]=-\lambda Y$ and ${\rm div}(Y)=\big(A(\beta_1+1)+B(\beta_2+1)\big)x^{\beta_1}y^{\beta_2}$, so the equation becomes
$$\Big(-\lambda\,{\rm Id}+2\binom{\lambda_1}{\lambda_2}(\beta_1+1,\ \beta_2+1)\Big)\binom{A}{B}=0.$$
If both monomials exist, the determinant of this matrix is $-\lambda\big(\lambda_1(\beta_1+2)+\lambda_2(\beta_2+2)\big)$. As $\lambda_1,\lambda_2$ are linearly independent over $\mathbb{F}_p$, it vanishes only for $\beta\equiv(0,0)$ and $\beta\equiv(-2,-2)$ modulo $p$, and in these two cases the matrix has rank one. If only one monomial exists, then either $\beta_1=-1$ or $\beta_2=p-1$ (only the first monomial), or $\beta_2=-1$ or $\beta_1=p-1$ (only the second), and the corresponding scalar equals $-\lambda$, $\lambda_1(\beta_1+2)-\lambda_2\beta_2$, $-\lambda$ or $\lambda_2(\beta_2+2)-\lambda_1\beta_1$, respectively; each is a nonzero $\mathbb{F}_p$-combination of $\lambda_1,\lambda_2$, hence nonzero. Hence $\L_{J^{-1}X_0}$ is two-dimensional. Applying the grading torus and multiplying by $t^{2p-4}$ we obtain a curve in the $\k^\times G$-orbit of $\chi$ whose limit at $t=0$ is $J^{-1}X_0$, and since $\L^*_{\rm sing}$ is a closed cone, $\chi$ is regular as well.

Now suppose $\lambda_1,\lambda_2$ are dependent. We show $(d\psi)_X(Y)=0$ for all $Y\in\L$, in three steps.

{\it Step 1.} For $A\in\gl(2,\k)=\L_0$ let $\varphi(A)$ be the coefficient of $T^p$ in $\det(T\,{\rm Id}-A)$, where $A$ acts on $\O$. Since $X$ preserves the $\m$-adic filtration and acts on $\m^k/\m^{k+1}$ as $X_0$, we have $\psi(X')=\varphi(X'_0)$ for every $X'\in\L_{(0)}$ with linear part $X'_0$. Hence $(d\psi)_X(Y)=0$ for $Y\in\L_{\ge1}$, and $(d\psi)_X(Y)=(d\varphi)_{X_0}(Y)$ for $Y\in\L_0$.

{\it Step 2: $(d\varphi)_A=0$ whenever the eigenvalues of $A$ are dependent over $\mathbb{F}_p$.} For $A={\rm diag}(\alpha,\beta)$ the eigenvalues of $A$ on $\O$ are $i\alpha+j\beta$, $i,j\in\mathbb{F}_p$; let $\varphi_0(\alpha,\beta)$ be the coefficient of $T^p$ in $\prod_{i,j}(T-i\alpha-j\beta)$. The multiset $\{i\alpha+j\beta\}$ is invariant under $(\alpha,\beta)\mapsto(\beta,\alpha)$ and under $(\alpha,\beta)\mapsto(\alpha,\beta+t\alpha)$, $t\in\mathbb{F}_p$, and so is $\varphi_0$. At $(\alpha,0)$ we have $\prod_{i,j}(T-i\alpha)=(T^p-\alpha^{p-1}T)^p$, so $\varphi_0(\alpha,0)=-\alpha^{p(p-1)}$ and $\del_\alpha \varphi_0(\alpha,0)=0$; and $\del_\beta \varphi_0(\alpha,0)=0$ because the derivative in $\beta$ of $\prod_{i,j}(T-i\alpha-j\beta)$ at $\beta=0$ is $-\sum_{i,j}j\prod_{(i',j')\ne(i,j)}(T-i'\alpha)$, and $\sum_jj=0$. Differentiating the invariance under $(\alpha,\beta)\mapsto(\alpha,\beta+t\alpha)$ transports this to the points $(\alpha,t\alpha)$, and the symmetry gives the points $(0,\beta)$. Thus the gradient of $\varphi_0$ vanishes at every pair of dependent eigenvalues. Since $\varphi_0$ is symmetric, $\varphi_0(\alpha,\beta)=\Phi(\alpha+\beta,\alpha\beta)$ for a polynomial $\Phi$, and $\varphi(A)=\Phi({\rm tr}\,A,\det A)$ for all $A$, both sides being conjugation invariant and equal on diagonal matrices. Hence $(d\varphi)_A=\Phi_1\,d\,{\rm tr}+\Phi_2\,d\det$, with the partial derivatives $\Phi_1,\Phi_2$ evaluated at $({\rm tr}\,A,\det A)$, and it suffices to show $\Phi_1=\Phi_2=0$ at $(\alpha+\beta,\alpha\beta)$ for dependent $\alpha,\beta$. As $\del_\alpha \varphi_0=\Phi_1+\beta \Phi_2$ and $\del_\beta \varphi_0=\Phi_1+\alpha \Phi_2$, this follows from the vanishing of the gradient if $\alpha\ne\beta$. If $\alpha=\beta=\lambda$, it remains to show $\Phi_2(2\lambda,\lambda^2)=0$. Using the invariance with $t=-1$,
$$\Phi(2\lambda,\lambda^2-\delta^2)=\varphi_0(\lambda+\delta,\lambda-\delta)=\varphi_0(\lambda+\delta,-2\delta),$$
and the coefficient of $\delta^2$ on the right is $\frac12\big(\del_\alpha^2\varphi_0-4\del_\alpha\del_\beta \varphi_0+4\del_\beta^2\varphi_0\big)(\lambda,0)$. Here $\del_\alpha^2\varphi_0(\lambda,0)=\del_\alpha\del_\beta \varphi_0(\lambda,0)=0$ by the formulas above, and $\del_\beta^2\varphi_0(\lambda,0)=0$ because the second derivative in $\beta$ at $\beta=0$ of $\prod_{i,j}(T-i\alpha-j\beta)$ is $P_0\big((\sum_{i,j}\frac{j}{T-i\alpha})^2-\sum_{i,j}\frac{j^2}{(T-i\alpha)^2}\big)$, $P_0=\prod_{i,j}(T-i\alpha)$, and $\sum_jj=\sum_jj^2=0$. Hence $\Phi_2(2\lambda,\lambda^2)=0$.

{\it Step 3: $Y\in\L_{-1}$.} Conjugation invariance gives $(d\psi)_X([D,X])=0$ for $D\in\L$. If $D\in\L_{-1}$, then $[D,X]=[D,X_0]+[D,X_1]+Y'$ with $[D,X_1]\in\L_0$ and $Y'\in\L_{\ge1}$, so Steps 1 and 2 give $(d\psi)_X([D,X_0])=0$. This settles all $Y$ in the image of $\ad X_0$ on $\L_{-1}$, which is $\L_{-1}$ unless an eigenvalue of $X_0$ vanishes. If $X_0=0$, then $(d\psi)_X(Y)=-{\rm tr}(X^{p^2-p-1}Y)$ involves at least one factor $X_0$ in each term of degree zero, and vanishes. If $X_0=\lambda_2y\del_y$ with $\lambda_2\ne0$ (the case $X_0=\lambda_1x\del_x$ is symmetric), it remains to treat $Y=\del_x$. By \eqref{eq:2.3}, $(d\psi)_X(\del_x)=-{\rm tr}(X^{p^2-p-1}\del_x)$, and only the terms $X_0^rX_1X_0^{p^2-p-2-r}\del_x$ contribute. If $X_1$ maps $x^{i-1}y^j$ to $(\alpha(i-1)+\beta j)x^iy^j$ modulo other monomials, the trace equals $-\sum_{i,j}i\big(\alpha(i-1)+\beta j\big)(p^2-p-1)(j\lambda_2)^{p^2-p-2}$, which vanishes since $\sum_ii=\sum_ii(i-1)=0$ for $p\ge5$. If $X_0$ is nilpotent and nonzero, every term $X_0^rX_1X_0^{p^2-p-2-r}$ contains a power of $X_0$ of exponent at least $p$, and $X_0^p=0$ on $\O$. Thus $(dF_\psi)_\chi=0$. Since $dF_1\wedge dF_2=dF_1\wedge dF_\psi=\gamma\Pi$ by Lemma~\ref{L:2.5}, all submaximal Pfaffians of the matrix $\big(\chi([e_i,e_j])\big)$ vanish. Hence this skew-symmetric matrix has rank at most $n-4$, so that $\dim\L_\chi\ge4$ and $\chi$ is singular.
\end{proof}
If $\chi\in\L^*_{\rm reg}$ and $J\chi\in\L_{(0)}$, then Lemma~\ref{L:4.3} shows that $t^{2p-4}\alpha_t\cdot\chi\to J^{-1}X_0$ for $\alpha_t(x)=tx$, $\alpha_t(y)=ty$, and the limit is regular. Lemma~\ref{L:4.1} again yields $\chi\notin\mathcal{X}_{d_0}$. Combining this with \ref{ss:4.2} we obtain $\mathcal{X}_{\rm sing}\subseteq\L^*_{\rm sing}$, so that $\Psi^{-1}(\L^*_{\rm reg})\subseteq\Omega$.

\subsection{}
\label{ss:4.4}
It remains to show that no point of $\Psi^{-1}(\L^*_{\rm sing})$ lies in $\Omega$.
\begin{Lemma}\label{L:4.4}
Let $X=a\del_x+b\del_y\in\L$ with $b\in y^2\O$. Then ${\rm tr}_\O(X^\ell Y)=0$ for all $Y\in\L$ and all $\ell\ge 0$. In particular, $(dF_\psi)_{J^{-1}X}=0$.
\end{Lemma}
\begin{proof}
Put $R=\k[x]/(x^p)$ and decompose $\O=\bigoplus_{j=0}^{p-1}Ry^j$; we say that a $\k$-linear operator on $\O$ has degree $d$ if it maps $Ry^j$ into $Ry^{j+d}$ for all $j$. Thus $f\del_x$ with $f\in R$ has degree $0$ and $\del_y$ has degree $-1$. Write $a=a_0+a_1y+\cdots$, $b=b_2y^2+\cdots$. On each layer the degree-zero block of $X$ is $a_0\del_x$, and its degree-one block from layer $j$ to layer $j+1$ is $a_1\del_x+jb_2$. Only the components of $Y$ of degrees $0$ and $-1$ contribute to the trace. For $Y=h\del_x+yk\del_y$ the trace equals $\sum_j{\rm tr}_R\big((a_0\del_x)^\ell(h\del_x+jk)\big)=0$. For $Y=h\del_y$ the diagonal block on layer $j$ is $j(S_0+(j-1)S_1)h$ with $S_0,S_1$ independent of $j$, and its total trace vanishes since $\sum_jj=\sum_jj(j-1)=0$ for $p\ge 5$. As the adjugate of $T\,{\rm Id}-X$ is a polynomial in $X$, the last statement follows.
\end{proof}
\begin{Proposition}\label{P:4.5}
If $\chi\in\L^*_{\rm sing}$ then $(dF_2)_\chi=0$ and $w(\chi)=0$.
\end{Proposition}
\begin{proof}
Since $d(h_0^p)=0$, we have $dF_2=dF_\psi$. By Lemma~\ref{L:3.7}, $F_2$ and $w$ are $G$-semi-invariants, so the zero sets of $dF_2$ and $w$ are $G$-stable. If $J\chi\notin\L_{(0)}$, Lemma~\ref{L:4.2} allows us to assume that $J\chi$ has the form \eqref{eq:4.2} with $c\in y^2\k[y]$. Then $J\chi(y)\in y^2\O$, and Lemma~\ref{L:4.4} gives $(dF_2)_\chi=0$; as both $J\chi(y)$ and $(J\chi)^{[p]}(y)$ lie in $y^2\O$, we also have $w(\chi)=0$. If $J\chi\in\L_{(0)}$, then Lemma~\ref{L:4.3} gives $(dF_2)_\chi=0$, and $w(\chi)=0$ because $J\chi$ and $(J\chi)^{[p]}$ vanish at the origin.
\end{proof}
\begin{proof}[Proof of Theorem~\ref{main3}]
We have shown that $\Psi^{-1}(\L^*_{\rm reg})\subseteq\Omega$. Identify ${\rm Specm}\,Z$ with the subvariety of $\A^{n+2}$ given by \eqref{eq:1}, with coordinates $(z_b)_b,u,v$. Let $P_0\in{\rm Specm}\,Z$ have coordinates $z_b=\chi(e_b)^p$, $u=\alpha$, $v=\beta$, where $\chi\in\L^*_{\rm sing}$, and put $f_1=u^p+uv-F_1(\underline{z})$ and $f_2=v^p-w(\underline{z})^{p-2}u-F_2(\underline{z})$. Since $F_2$ and $w$ have coefficients in $\mathbb{F}_p$, we have $\big(\del F_2(\underline{z})/\del z_b\big)(P_0)=\big((\del F_2/\del e_b)(\chi)\big)^p$, and $w(\underline{z})$ takes the value $w(\chi)^p$ at $P_0$. Proposition~\ref{P:4.5} now shows that
$$(df_2)_{P_0}\,=\,\big(-dF_2(\underline{z})-(p-2)w(\underline{z})^{p-3}u\,dw(\underline{z})-w(\underline{z})^{p-2}du\big)_{P_0}\,=\,0.$$
As ${\rm Specm}\,Z$ is a complete intersection of codimension two by Theorem~\ref{main2}, its Jacobian matrix has rank $<2$ at $P_0$. Thus $P_0$ is a singular point of ${\rm Specm}\,Z$, and hence $P_0\notin\Omega$ by \cite[6.4]{PS99}.
\end{proof}
In particular, $U_\chi(\L)$ has a simple module of dimension $<N$ for every $\chi\in\L^*_{\rm sing}$, so that $\mathcal{X}_{\rm sing}=\L^*_{\rm sing}$.
Combining Theorem~\ref{main3} with Lemmas~\ref{L:4.2} and \ref{L:4.3} we obtain an explicit description of the image of the Azumaya locus.
\begin{Corollary}\label{C:4.6}
A linear function $\chi\in\L^*$ lies in $\Psi(\Omega)$ if and only if one of the following holds:

\smallskip

\noindent
{\rm (a)} $J\chi\notin\L_{(0)}$ and in a normal form \eqref{eq:4.2} of $\chi$ one has $c(0)\ne 0$ or $c'(0)\ne 0$;

\smallskip

\noindent
{\rm (b)} $J\chi\in\L_{(0)}$ and the eigenvalues of the linear part of $J\chi$ are linearly independent over $\mathbb{F}_p$.

\smallskip

\noindent
Moreover, $\Psi^{-1}(\chi)\subseteq\Omega$ for every $\chi\in\Psi(\Omega)$, the complement $\L^*\setminus\Psi(\Omega)$ has codimension $\ge 2$ in $\L^*$, and $(dF_2)_\chi=0$ and $w(\chi)=0$ for all $\chi\notin\Psi(\Omega)$.
\end{Corollary}
\begin{proof}
By Theorem~\ref{main3}, $\Omega=\Psi^{-1}(\L^*_{\rm reg})$, so that $\Psi(\Omega)=\L^*_{\rm reg}$ because $\Psi$ is surjective. Lemmas~\ref{L:4.2} and \ref{L:4.3} identify $\L^*_{\rm reg}$ with the set described in (a) and (b). The codimension statement was observed after Lemma~\ref{L:4.2}, and the last statement is Proposition~\ref{P:4.5}.
\end{proof}
In particular, every $\chi$ whose stabiliser $\L_\chi$ is a torus lies in $\Psi(\Omega)$, and so does the $T_0$-unstable linear function $\eta$.

\subsection{}
\label{ss:4.5}
We now determine the geometric monodromy group of $\Psi$, that is, the Galois group of the normal closure of the field extension $Q(\L)/Q_p(\L)$. The Jacobian determinant of the two relations \eqref{eq:1} with respect to $(u,v)$ is $\det\begin{pmatrix}v&u\\-w(\underline{z})^{p-2}&0\end{pmatrix}=uw(\underline{z})^{p-2}$, so $\Psi$ is \'etale where $uw(\underline{z})\ne0$. It is invertible on the open set where $F_1(\underline{z})w(\underline{z})$ does not vanish, because
$$u\big(u^{p-1}+v\big)\,=\,F_1(\underline{z}).$$ On the slice \eqref{eq:3.5} put $\alpha=a^p$ and $\beta=b^p$. The equations of $Z$ become
$$\alpha=u^p+uv,\qquad\beta=v^p-u,$$
and setting $\beta=0$ we obtain the connected cover of degree $p^2$ (whose monodromy group is a subgroup of that of $\Psi$)
$$v^{p^2}+v^{p+1}\,=\,\alpha.$$
At $\alpha=0$ it has a tamely ramified branch of index $p+1$ at $v=0$ and $p^2-p-1$ unramified branches; hence inertia is generated by a $(p+1)$-cycle. 
The monodromy group is primitive. Indeed, a nontrivial block system would consist of $p$ blocks of size $p$. A permutation preserving such a system and moving two blocks moves at least $2p$ letters, whereas a cycle preserving each block has length at most $p$; neither is compatible with a $(p+1)$-cycle. 

By Jones's classification of primitive groups containing a cycle \cite[Corollary~1.3]{Jones}, the group contains the alternating group $\mathfrak{A}_{p^2}$, and as $p+1$ is even, the $(p+1)$-cycle is odd. Therefore the monodromy group is $\mathfrak{S}_{p^2}$. Thus the normal closure $F$ of $Q(\L)/Q_p(\L)$ has Galois group $\mathfrak{S}_{p^2}$, acting on the $p^2$ embeddings of $Q(\L)$ into $F$ over $Q_p(\L)$. Hence $Q(\L)=F^H$, where $H\cong\mathfrak{S}_{p^2-1}$ is a point stabiliser. As $H$ is a maximal subgroup of $\mathfrak{S}_{p^2}$, the extension $Q(\L)/Q_p(\L)$ has no proper intermediate fields (by the Galois correspondence). This agrees with Kadets's classification of generic trinomial monodromy \cite[Theorem~5.4]{Kadets}.

\subsection{}
\label{ss:4.6}
Finally, we explain why Theorem~\ref{main3} fails in characteristics $2$ and $3$. If $\Char\k=2$ then $\L\cong\sl(3)$, so that $\L$ is the Lie algebra of a simply connected group. As explained in \ref{ss:1.3}, in this case $\Psi(\Omega)=\L^*$ by \cite[Sect.~4]{FP90}, whereas $\L^*_{\rm sing}\ne\varnothing$.

Now let $\Char\k=3$. The definitions of $F_1,F_\psi,h_0,w$ still make sense, and the residue computations of Section~\ref{s:2} show that $F_1,F_2\in S(\L)^\L$ and that $(dF_1)_\eta$ and $(dF_2)_\eta$ are linearly independent. The elimination in the proof of Lemma~\ref{L:4.2} only involves the denominators $1$ and $2$, and a nonzero constant or linear coefficient of $c$ still gives a two-dimensional stabiliser. Hence $\L^*_{\rm sing}$ has codimension $\ge 2$, and the proof of Theorem~\ref{T:2.6} yields $S(\L)^\L=S(\L)^3[F_1,F_2]$ with $\deg F_1=4$ and $\deg F_2=6$.

The arguments of Section~\ref{s:3} go through for $p=3$ with two modifications.

\smallskip

\noindent
(i) Lemma~\ref{L:3.1}(i) holds with $2p-4=2$, and this suffices for Lemmas~\ref{L:3.1}(ii), \ref{L:3.2} and Proposition~\ref{P:3.3}, which only involve derivatives of $a_1,a_2$ of total order at most $2$. Hence $u=\R(\D^3)$ is central with ${\rm gr}(u)=F_1$, and the symmetrisation in \ref{ss:3.5} applies without change.

\smallskip

\noindent
(ii) In \ref{ss:3.3}, Lemma~\ref{L:3.4} still gives $[Y,v']\in U_4$ for all $Y\in\L$, but the final degree estimate fails by one, since $p^2-4p+2=-1$. Put $E_Y=[Y,v']$. As $v'$ has integer degree $12$ and every element of $U_4$ has integer degree $\le 12$, we have $E_Y=0$ for $Y\in\L_{\ge 1}$. Since $T_0$ acts on $v'$ by the character $\det^{6}$, we also have $E_Y=0$ for $Y\in\t_0$. Let $\del\in\L_{-1}$. The subalgebra $\L_{\ge 2}=\L_2\oplus\L_3$ is abelian, and comparing PBW and integer degrees we see that $E_\del$ lies in $S^3(\L_3)\L_2\subset U(\L_{\ge2})$; it is a $T_0$-weight vector of weight $(5,6)$ for $\del=\del_x$ and $(6,5)$ for $\del=\del_y$. If $Z\in\L_1$ and $[Z,\del]\in\t_0$, then $[Z,E_\del]=E_{[Z,\del]}+[\del,E_Z]=0$. Put $A=x^2y^2\del_x$ and $B=x^2y^2\del_y$. For $\del=\del_x$ the weight space is spanned by $A^3\cdot x^2y\del_y$, $A^2B\cdot x^2y\del_x$, $A^2B\cdot xy^2\del_y$ and $AB^2\cdot xy^2\del_x$; write $E_{\del_x}=e_1A^3\cdot x^2y\del_y+e_2A^2B\cdot x^2y\del_x+e_3A^2B\cdot xy^2\del_y+e_4AB^2\cdot xy^2\del_x$. The elements $Z=x^2\del_x,\,y^2\del_x,\,xy\del_y,\,y^2\del_y$ satisfy $[Z,\del_x]\in\t_0\cup\{0\}$, and $\ad Z$ kills $\L_3$ and acts on $\L_2$ by
$$\begin{array}{l}
[x^2\del_x,xy^2\del_x]=-A,\quad [x^2\del_x,xy^2\del_y]=\ \ \ B,\quad [y^2\del_x,x^2y\del_y]=A,\\[2pt]
[xy\del_y,x^2y\del_x]=-B,\quad [xy\del_y,xy^2\del_x]=\,\, 2A,\quad [xy\del_y,xy^2\del_y]=B,\\[2pt]
[y^2\del_y,x^2y\del_x]=\ \ \ A,\quad [y^2\del_y,x^2y\del_y]=-B,
\end{array}$$
all other brackets of these $Z$ with the four elements of $\L_2$ being zero. The equations $[Z,E_{\del_x}]=0$ therefore read $e_3-e_4=0$, $e_1=0$, $-e_2+e_3+2e_4=0$ and $-e_1+e_2=0$. Their solutions are $e_1=e_2=0$, $e_3=e_4$, that is, multiples of
$$A^2B\cdot xy^2\del_y+AB^2\cdot xy^2\del_x\,=\,[\del_x,A^2B^2].$$
By symmetry, $E_{\del_y}$ is a multiple of $[\del_y,A^2B^2]$. Write $E_{\del_x}=\lambda[\del_x,A^2B^2]$ and $E_{\del_y}=\lambda'[\del_y,A^2B^2]$. Since $[\del_x,E_{\del_y}]=[\del_y,E_{\del_x}]$ and the principal symbol of $[\del_x,[\del_y,A^2B^2]]$ is nonzero, we get $\lambda=\lambda'$. Then $v''=v'-\lambda A^2B^2$ commutes with $\L_{-1}$, with $\L_{\ge1}$ (as $[\L_1,\L_3]=0$) and with $\L_0=[\L_1,\L_{-1}]$. Hence $v=v''+h_0(\underline{z})$ is central with ${\rm gr}(v)=F_2$.

\smallskip

\noindent
The rest of Section~\ref{s:3} applies verbatim, except that in the proof of Proposition~\ref{P:3.6} the degree count admits one additional term $\gamma\,z_Az_B$ in the expansion of $u^3-F_1(\underline{z})$. Its coefficient is a $G$-semi-invariant of degree $2$ and weight $3(p-2)$, and \eqref{eq:3.6} has no solution with $q=3$, so that $\gamma=0$. We obtain
$$Z\,=\,\textstyle{\bigoplus_{0\le i,j\le 2}}\,Z_p(\L)u^iv^j,\qquad u^3+uv\,=\,F_1(\underline{z}),\qquad v^3-w(\underline{z})u\,=\,F_2(\underline{z}).$$
The Azumaya and smooth loci of $Z$ coincide by \cite[Theorem~3.13]{BM}. Indeed, $U$ has finite global dimension, it is a maximal Cohen--Macaulay $Z$-module, and it is Azumaya in codimension one by \cite[Proposition~5.2]{PS99}. Its PI-degree is $3^8$ by \cite[Theorems~5 and 6]{Zass}.

Take $J\chi=(1+x^2)\del_x$. Then
$$\L_\chi\,=\,\{(1-x^2)c(y)\del_x+b\del_y\,|\,\,c\in\k[y]/(y^3),\ b\in\k\}$$
has dimension $4$, so $\chi\in\L^*_{\rm sing}$. Put $X=J\chi$. Since $X(x)=1+x^2$, $X^2(x)=2x$, $X^3(x)=2(1+x^2)=-X(x)$ and $X(y)=0$, we have $X^{[3]}=-X$. As ${\rm div}(X)=2x$, $q_X=(X+x)^2(1)=X(x)+x^2=1-x^2$. On each of the three summands $\k[x]/(x^3)\,y^j$ the operator $X$ has matrix with characteristic polynomial $T^3+T$ in the basis $y^j,xy^j,x^2y^j$, so $\det(T\,{\rm Id}-X)=T^9+T^3$. Thus $F_1(\chi)=0$, $F_2(\chi)=F_\psi(\chi)+h_0(\chi)^3=2$ and $w(\chi)=0$. The three points of $\Psi^{-1}(\chi)$ are $(u,v)=(0,2),(1,2),(-1,2)$. Moreover, $(dF_2)_\chi=0$: by \eqref{eq:2.3} we have $(d\psi)_X(Y)=-{\rm tr}(X^5Y)=-{\rm tr}(XY)$ for $X=J\chi$, and this vanishes after summing over the three $y$-layers of $\O$. On the other hand, for $Y=-x^2\del_y$ the derivative of $X\mapsto X^3$ in the direction $Y$ sends $y$ to $X^2(-x^2)=1+x^2$, so that $(dw)_\chi(Y)=1$. At the points with $u=\pm 1$ the Jacobian matrix of the two relations has rank two, so these points are smooth and therefore lie in $\Omega$, although $\chi$ is singular. The point with $u=0$ is singular.

\end{document}